\documentclass[12pt]{amsart}
\usepackage{amscd}
\let\oldlabel=\label
\def\prellabel{\marginparsep=1em\marginparwidth=44pt
    \def\label##1{\oldlabel{##1}\ifmmode\else\ifinner\else
         \marginpar{{\footnotesize\ \\ \tt
                    ##1}}\fi\fi}}
\advance\textheight1cm
\def\Ker{{\operatorname{Ker}}}
\def\det{{\operatorname{det}}}
\def\NN{{\mathbb N}}

\def\Cl{{\operatorname{Cl}}}

\def\ggr.aut{{\operatorname{gr.aut}}}

\def\Gr.Hom{{\operatorname{Gr.Hom}}}
\def\L{\operatorname{L}}
\def\width{\operatorname{width}}
\def\int{{\operatorname{int}}}
\def\conv{\operatorname{conv}}
\def\het{\operatorname{ht}}
\def\E{\operatorname{E}}
\def\N{\operatorname{N}}
\def\Spec{\operatorname{Spec}}
\def\gr{\operatorname{gr}}
\def\gp{\operatorname{gp}}
\def\rank{\operatorname{rank}}
\def\Ker{{\operatorname{Ker}}}
\def\codim{\operatorname{codim}}
\def\Proj{\operatorname{Proj}}
\def\maxSpec{\operatorname{maxSpec}}
\def\Im{\operatorname{Im}}
\def\Pic{\operatorname{Pic}}

\def\Vect{\operatorname{Vect}}
\def\Pol{\operatorname{Pol}}
\def\GL{\operatorname{GL}}
\def\join{\operatorname{join}}

\def\Col{\operatorname{Col}}
\let\tensor=\otimes
\def\RR{{\mathbb R}}
\def\CC{{\mathbb C}}
\def\QQ{{\mathbb Q}}
\def\ZZ{{\mathbb Z}}
\def\NN{{\mathbb N}}
\def\TT{{\mathbb T}}

\def\EE{{\mathbb E}}
\def\AA{{\mathbb A}}
\def\PP{{\mathbb P}}
\def\Oo{{\mathcal O}}

\def\pp{{\mathfrak p}}

\let\epsilon=\varepsilon
\let\phi=\varphi
\let\theta=\vartheta

\makeatletter
\def\>#1>#2>{\setbox0\hbox{$\scriptstyle
 \;{#1}\;\;$}\setbox1\hbox{$\scriptstyle\;{#2}\;\;$}\setbox2
 \hbox{$#2$}
 \global\bigaw@\minaw@
 \ifdim\wd0>\bigaw@\global\bigaw@\wd0\fi
 \ifdim\wd1>\bigaw@\global\bigaw@\wd1\fi
 \ifdim\wd2>\z@
 \mathrel{\mathop{\hbox to\bigaw@{\rightarrowfill}}\limits^{#1}_{#2}}\else
 \mathrel{\mathop{\hbox to\bigaw@{\rightarrowfill}}\limits^{#1}}\fi
}
\makeatother

\def\vertex{\circle*{0.15}}

\advance\headheight1.15pt

\newtheorem{lemma}{Lemma}[section]

\newtheorem{theorem}[lemma]{Theorem}
\newtheorem{proposition}[lemma]{Proposition}

\theoremstyle{definition}

\newtheorem{remark}[lemma]{Remark}
\newtheorem{example}[lemma]{Example}

\begin{document}

\title[Polytopal linear retractions]{Polytopal linear retractions}

\author{Winfried Bruns \and Joseph Gubeladze}
\address{Universit\"at Osnabr\"uck,
FB Mathematik/Informatik, 49069 Osnabr\"uck, Germany}
\email{Winfried.Bruns@mathematik.uni-osnabrueck.de}
\address{A. Razmadze Mathematical Institute, Alexidze St. 1, 380093
Tbilisi, Georgia}
\email{gubel@rmi.acnet.ge}
\thanks{The second author was supported by the Max-Planck-Institut f\"ur Mathematik
in Bonn and INTAS}
\subjclass{Primary 13F20, 14M25; Secondary 52C07}
\keywords{Polytopal algebra, retracts, affine semigroup ring, binomial ideal}

\begin{abstract}

We investigate graded retracts of polytopal algebras (essentially the
homogeneous rings of affine cones over projective toric varieties) as polytopal
analogues of vector spaces. In many cases we show that these retracts are again
polytopal algebras and that codimension $1$ retractions factor through
retractions preserving the semigroup structure. We expect that these results
hold in general.

This paper is a part of the project started in \cite{BG1,BG2}, where we have
investigated the graded automorphism groups of polytopal algebras. Part of the
motivation comes from the observation that there is a reasonable `polytopal'
generalization of linear algebra (and, subsequently, that of algebraic
$K$-theory).
\end{abstract}

\maketitle

\section{Introduction}

The category $\Vect_k$ of finite-dimensional vector spaces over a field $k$ has
a natural extension that we call the {\em polytopal $k$-linear category}
$\Pol_k$. The objects of $\Pol_k$ are the {\em polytopal semigroup algebras} (or
just {\em polytopal algebras}) in the sense of \cite{BGT}. That is, an object
$A\in|\Pol_k|$ is (up to graded isomorphism) a standard graded $k$-algebra
$k[P]$ associated with a lattice polytope $P\subset\RR^n$. The morphisms of
$\Pol_k$ are the homogeneous $k$-algebra homomorphisms. Here, as usual,
`standard graded' means `graded and generated in degree 1', and a lattice
polytope is a finite {\em convex} polytope whose vertices are lattice points of
the integral lattice $\ZZ^n\subset\RR^n$.

The generators of $k[P]$  correspond bijectively to the lattice points in $P$,
and their relations are the binomials representing the affine dependencies of the
lattice points. More precisely, let $\L_P$
denote the set of lattice points in $P$. Then $S_P$ is the sub-semigroup
of $\ZZ^{n+1}$ generated by
$$
\E_P=\{(x,1)|\ x\in\L_P\}.
$$
and $k[P]=k[S_P]$ is the semigroup algebra of $S_P$ with coefficients in $k$ (in
[BGT, BG1, BG2] we have always used the notation $k[S_P]$). Unless specified
otherwise, we will identify the sets $\L_P$ and $\E_P$.

The category $\Pol_k$ contains $\Vect_k$ as a full subcategory. In fact, we can
identify a vector space with the degree $1$ component of its symmetric algebra,
which, upon the choice of a basis, can be considered as a polynomial ring
$k[X_1,\dots,X_n]$. This polynomial ring is isomorphic to the polytopal algebra
$k[{\Delta_{n-1}}]$ defined by the $(n-1)$-simplex $\Delta_{n-1}$. Vector
space homomorphisms extend to homomorphisms of symmetric algebras, and thus to
homomorphisms of the corresponding polytopal algebras.

Observe that polytopal algebras are essentially the homogeneous
coordinate rings of projective normal toric varieties. Actually,
one needs the extra condition of {\em very ampleness} for the underlying
polytope to get such a homogeneous coordinate ring (see [BG1,\S5] for details).

Our investigation of polytopal algebras is motivated by two closely related
goals: (1) to find the connections between the combinatorial structure of $P$
and the algebraic structure of $k[P]$, and (2) to extend theorems valid in
$\Vect_k$ to $\Pol_k$.

It follows from \cite{Gu} that an algebra isomorphism of $k[P]$ and $k[Q]$
implies the isomorphism of $P$ and $Q$ as lattice polytopes. This result
identifies the objects of the category $\Pol$ of lattice polytopes with the
objects of $\Pol_k$ (up to isomorphism), but there remains the question to which
extent the morphisms in $\Pol_k$ are determined by those in $\Pol$, namely the
$\ZZ$-affine maps between lattice polytopes. Similarly one must ask whether
certain classes of morphisms in $\Pol_k$ can be described in the same way as the
corresponding classes in $\Vect_k$.

As we have shown in \cite{BG1}, there is a total analogy with the linear
situation for the automorphism groups in $\Pol_k$ (called {\em polytopal linear
groups} in \cite{BG1}): they are generated by elementary automorphisms (generalizing
elementary matrices), toric automorphisms (generalizing diagonal invertible
matrices) and automorphisms of the underlying polytope; moreover, there are
normal forms for such representations of arbitrary automorphisms. (In \cite{BG2} this
analogy has been extended to automorphisms of so-called polyhedral algebras.)

In the present paper we study {\em retractions} in $\Pol_k$, that is, idempotent
homogeneous endomorphisms of polytopal algebras and their images. Our results
support the following two conjectures:
\medskip

\noindent {\bf Conjecture A.} {\em The images of retractions in $\Pol_k$ are
polytopal algebras.}
\medskip

\noindent {\bf Conjecture B.} {\em A codimension 1 retraction factors through either a
facet projection or an affine lattice retraction of the underlying
lattice polytope.}
\medskip

These conjectures generalize the standard facts that every finitely generated
vector space has a basis and that an idempotent matrix is conjugate to a
sub-unit matrix (ones and zeros on the main diagonal and zeros anywhere else).
Conjecture B must be restricted to codimension 1 since there exist
counterexamples for higher codimension.

In the following we say that a retraction $h$ of a polytopal algebra is {\em
polytopal} if $\Im(h)$ is a polytopal algebra, and codimension $1$ retractions
satisfying the stronger condition of Conjecture B are called {\em tame} (the
notion of tameness has a natural generalization to higher codimension). We will
provide evidence supporting Conjectures A and B by proving them in number of
special cases (we always assume that the field $k$ is algebraically closed):
\begin{itemize}
\item[(1)]
If $\dim \Im(h)\le 2$, then $h$ is polytopal (Theorem \ref{2.2}).
\item[(2)]
Retractions of normal polytopal algebras $C$ that arise from a tensor product
decomposition $A\otimes B=C$ of graded $k$-algebras are polytopal; in other
words, graded tensor factors of normal polytopal algebras are again polytopal
(Theorem \ref{3.1}).
\item[(3)]
If $\Im(h)=k[Q]$ for a lattice subpolytope $Q$ of $P$ such that $Q$ meets the interior of $P$, then $h$ is tame (Theorem \ref{7.2a}).
\item[(4)]
If $\Im(h)=k[Q]$, but the subpolytope $Q$ is contained in the boundary of $P$, then $h$ need not be tame (Examples \ref{4.2} and \ref{WildTame}) in general. However, it is tame if $\codim(h)=1$ (Theorem \ref{7.2b}).
\item[(5)]
If $h$ is a codimension $1$ retraction such that $\ggr.aut_k(\Im(h))$ contains a torus of
dimension $\dim\Im(h)$ whose action can be extended algebraically to $k[P]$, then $h$ is tame
(Theorem \ref{7.5}).
\item[(6)]
If $\dim(P)\le 2$, then all codimension $1$ retractions of $k[P]$ are tame
(Theorem \ref{8.1}).
\end{itemize}

In Section \ref{seg} we discuss the class of {\em
segmentonomial} ideals, that is, ideals generated by polynomials $f$ whose
Newton polytope has dimension $\le 1$. Section \ref{aut} contains a review of
the results of \cite{BG1}.
\bigskip

\noindent {\em Conventions and terminology.}\enspace
The basic facts about affine semigroup rings can be found in \cite[Chap.\ 6]{BH} and \cite{Gi}.
For toric varieties we refer the reader to \cite{Fu,Oda}. Newton polytopes are treated in \cite{GKZ}.

By $k$ we will always denote a field, and we set $k^*=k\setminus\{0\}$. $\ZZ$, $\QQ$,
$\RR$ are the rational, integer and real numbers and $\ZZ_+$,
$\QQ_+$, $\RR_+$ are the sub-semigroups of non-negative elements.

An {\em affine} semigroup is a finitely generated subsemigroup of a free
abelian group. The group of differences of a semigroup $S$ will be
denoted by $\gp(S)$. An affine semigroup is {\em positive} if its
maximal subgroup is trivial. A positive affine semigroup can be embedded
into a free commutative semigroup $\ZZ^m_+$, thus the semigroup ring
$k[S]$ has a (non-unique) grading $k[S]=k\oplus A_1\oplus A_2\oplus\cdots$ making all
nontrivial elements of $S$ homogeneous of positive degree.
The uniquely determined minimal generating set of a positive affine
semigroup $S$ coincides with its set of indecomposable elements.

In $\RR\otimes\gp(S)$ the semigroup $S$ spans the cone $C(S)$.

An affine semigroup $S$ is called {\em homogeneous} if $S$ does not have
nontrivial invertibles and its irreducible elements are contained in some proper
hyperplane of $\RR\otimes\gp(S)$. This condition is equivalent to the
requirement that the semigroup ring $k[S]$ is a standard graded ring in which
$S$ consists of homogeneous elements. Alternatively, there is a lattice polytope
$P$ and a subset $X\subset\L_P$ such that $k[S]$ is isomorphic to the subalgebra
of $k[P]$ generated by $X$.

We put $\rank(S)=\rank(\gp(S))$. An affine semigroup $S$ is called
{\em normal} if
$$
\forall x\in\gp(S),\ (c\in\NN,\ cx\in S)\Rightarrow(x\in S).
$$
The normality of $S$ is equivalent to that of $k[S]$.

The elements of an affine semigroup $S$ will be called {\em terms} of
$k[S]$, while the elements of type $as\in k[S]$, $a\in k$, $s\in S$,
will be called {\em monomials}.

For an affine semigroup $S$ we have the embeddings
$$
S\subset\gp(S)\subset\RR\otimes\gp(S)\approx\RR^{\rank(S)}.
$$
Therefore, we can speak of the {\em convex hull} $\conv(X)$ of a subset
$X\subset S$ in $\RR\otimes\gp(S)$. For an element $f\in k[S]$ there is
a unique canonical representation
$$
f=\sum_{s\in S} a_s s,\qquad a_s\in k.
$$
The polytope
$$
\conv(\{s\in S\mid a_s\neq 0\})\subset\RR\otimes\gp(S)
$$
is called the {\em Newton polytope} of $f$, which we denote by $\N(f)$.
We will frequently use the formula
$$
\N(fg)=\N(f)+\N(g),
$$
where `+' denotes the {\em Minkowski sum} of convex polytopes:
$$
P+Q=\{x+y\mid x\in P,\ y\in Q\}\subset\RR^n,\ \ \ P\subset\RR^n,\
Q\subset\RR^n.
$$

As remarked already, a lattice polytope $P\subset\RR^n$ is a finite
convex polytope whose vertices belong to the integral lattice $\ZZ^n$.
It is called {\em normal} if $k[P]$ is normal.

Let $S\subset\ZZ^n$ be an affine semigroup. The semigroup operation in $S$ is 
written additively when we do linear algebra in $\ZZ^n$, and multiplicatively 
when $S$ is considered as the set of terms in $k[S]$.  

Finally, we adopt the blanket assumption that the field $k$ is {\em
algebraically closed}, unless specified otherwise.

\section{Retracts of dimension two}\label{two}

A {\em retract} of a $k$-algebra $A$ is an algebra $B$ such that
there exist $k$-homomorphisms $f:B\to A$ and $g:A\to B$ with $g\circ f=1_B$.
This is equivalent to saying that there is an endomorphism $h:A\to A$
such that $h^2=h$ and $\Im(h)\approx B$. We will call such $g$ and $h$
{\em retractions} and will frequently make passages between the two
equivalent definitions. Moreover, all the retractions considered below
are supposed to be graded. For a retraction $h$ as above we put
$$
\codim(h)=\dim(A)-\dim(B).
$$
The arguments used in the sequel need $k$ to be algebraically closed. We
do not know whether the following implication holds
for a retraction $h:A\to A$, $A\in|\Pol_k|$, over a not necessarily
algebraically closed field $k$:
$$
\bar k\otimes h \text{ is polytopal (tame)}\implies h \text{ is polytopal (tame)}
$$
($\bar k$ is the algebraic closure).

That polytopality is in general not an invariant property under
scalar extension/restriction is exhibited by the following

\begin{example}\label{2.1}
Consider the standard graded $\RR$-algebra
$$
A=\RR[X,Y,Z]/(X^2+Y^2+Z^2).
$$
Then $A$ is a factorial non-polytopal algebra over $\RR$ while $\CC\otimes A$ is
isomorphic to the polytopal algebra $\CC[{2\Delta_1}]$ defined by a lattice
segment $2\Delta_1$ of length $2$.

The factoriality of $A$ is proved in [Fo,\S11]. But the only
factorial polytopal algebras (over any field) are polynomial algebras --
an easy observation. Hence $A$ is not polytopal because it is singular at
the irrelevant maximal ideal. But we have the isomorphism
$$
\alpha:\CC[U^2,UV,V^2]=\CC[{2\Delta_1}]\to \CC\otimes A
$$
defined by $U^2\mapsto X+iY$, $V^2\mapsto X-iY$, $UV\mapsto iZ$.
\end{example}

Conjecture A  holds in Krull dimension $\le 2$:

\begin{theorem}\label{2.2}
A retract $B$ of a polytopal algebra $A$ is polytopal if $\dim B\le 2$.
\end{theorem}

The crucial step in the proof is

\begin{proposition}\label{2.3}
Let $k$ be an algebraically closed field and $A$ a standard graded $k$-algebra
of dimension $2$. If $A$ is a normal domain and the class group $\Cl(A)$ is
finitely generated, then $A$ is isomorphic to $k[{c\Delta_1}]$ as a graded
$k$-algebra for some $c\in\NN$ (as usual, $\Delta_1$ is the unit segment).
\end{proposition}

\begin{proof}
We have the projectively normal embedding of $\Proj(A)$ given by $A$.
Therefore, the projective curve $\Proj(A)$ is normal and thus smooth.
Consider the exact sequence
$$
0\to\ZZ\to \Cl(\Proj(A))\to \Cl(A)\to0.
$$
of Weil divisors arising from viewing
$\Spec(A)$ as a cone over $\Proj(A)$ (\cite[Ex. II.6.3(b)]{Ha}).
Since $\Cl(A)$ is finitely generated, so is $\Cl(\Proj(A))$. In
particular the Jacobian ${\mathcal J}(\Proj(A))\approx \Cl^0(\Proj(A))$ is
trivial ($\Cl^0$ denotes degree zero divisor classes). Therefore the
genus of $\Proj(A)$ is $0$, or equivalently $\Proj(A)\approx\PP^1_k$.
Using the normality of $A$ once again we get
$$
A\approx\bigoplus_{i=0}^\infty H^0(\PP^1_k,{\mathcal L}^{\otimes i})
$$
for some very ample line bundle $\mathcal L$ on $\PP^1_k$. But due to the
equality $\Pic(\PP^1_k)=\ZZ$
such a line bundle is a positive multiple of ${\Oo}(1)$, and hence $A$ is the
Veronese subalgebra of the polynomial algebra $k[{\Delta_1}]$ of some
level $c\in\NN$.
\end{proof}

\begin{proof}[Proof of Theorem \ref{2.2}]
In case $\dim(B)=1$ it is easy to see that $B\approx k[X]$.

Consider the case $\dim(B)=2$. We write $A=k[P]$ and denote the retraction $A\to B\subset A$ by $g$.
Consider the set
$(S_P\cap\Ker(g))\subset k[P]$ of monomials. There is a unique face $F\subset P$ such
that $(S_P\cap\Ker(g))=(S_P\setminus S_F)$ and $k[P]/(S_P\cap\Ker(g))$ is
naturally isomorphic to $k[F]$. Then $g$ is a composite of the two retractions
$$
k[P]\>\pi>>k[F]\>\rho>>B
$$
where $\rho$ is the homomorphism induced by $g$. Observe that $\rho$ is in fact
a retraction as it is split by $\pi|_B$.

Therefore we can from the beginning assume
that $(S_P\cap\Ker(g))=0$. In this situation $g$ extends (uniquely) to
the normalizations
\begin{equation}
\tag{1}\overline{k[P]}=k[\bar S_P]\>\bar g>>\bar B.
\end{equation}
This extension is given by
$$
\forall z\in\bar S_P,\ \ \bar g(z)=\frac{g(x)}{g(y)},\ \
\text{where}\ x,y\in S_P\ \text{and}\ z=\frac{x}{y}.
$$
It is known that the semigroup $S_{nP}$ is normal for all natural
numbers $n\geq\dim(P)-1$ \cite{BGT}. Therefore, by restricting the retraction
(1) to the $n$th Veronese subalgebra for such a number $n$, we get the
retraction
\begin{equation}
\tag{2}k[{nP}]\>{\bar g_n}>>\bar B_{(n)}.
\end{equation}
Let us show that $\Cl(\bar B_{(n)})$ is finitely generated for $n\geq\dim(P)-1$.
We choose a lattice point $x$ of $S_{nP}$ that is in the interior of the cone
$C(S_P)$. By localization (2) gives rise to the retraction
\begin{equation}
\tag{3}(x\bar g_n(x))^{-1}k[{nP}]\to(\bar g_n(x))^{-1}\bar B_{(n)}.
\end{equation}

Since $(x\bar g_n(x))^{-1}k[{nP}]$ is a localization of the
Laurent polynomial ring $x^{-1}k[{nP}]=k[\gp(S_{nP})]$, it
is a factorial ring. Then its retract $(\bar g_n(x))^{-1}\bar B_{(n)}$
is factorial as well (for example, see \cite{Cos}). By Nagata's theorem
\cite{Fo} $\Cl(\bar B_{(n)})$ is generated by the classes of the height 1
prime ideals of $\bar B_{(n)}$ containing $\bar g_n(x)$ -- a finite set.

It is also clear from (2) that $\bar B_{(n)}$ is generated in degree 1.
Consequently, by Proposition \ref{2.3} for each $n\geq\dim(P)-1$ there
is a natural number $c_n$ and an isomorphism
$$
\phi_n:\bar B_{(n)}\to k[{c_n\Delta_1}].
$$ We now fix such a number $n$. Restricting $\phi_n$ and $\phi_{n+1}$ to the
iterated Veronese subalgebra $\bar B_{(n(n+1))}=(\bar B_{(n)})_{(n+1)}=(\bar
B_{(n+1)})_{(n)}$ we obtain two isomorphisms of $\bar B_{(n(n+1))}$ with
$k[{c_{n(n+1)}\Delta_1}]$. It follows that there exists $c\in\NN$ with $c_n=cn$
and $c_{n+1}=c(n+1)$, and furthermore the restrictions of $\phi_n$ and
$\phi_{n+1}$ differ by an automorphism of $k[{c_{n(n+1)}\Delta_1}]$. However,
each automorphism of $k[{c_{n(n+1)}\Delta_1}]$ can be lifted to an automorphism
of $k[{\Delta_1}]$, and then restricted to all Veronese subrings of
$k[{\Delta_1}]$. (This follows from the main result of \cite{BG1}.) Therefore we can
assume that the restrictions of $\phi_n$ and $\phi_{n+1}$ coincide. Then they
define an isomorphism of the subalgebra $V$ of $\bar B$ generated by the
elements in degree $n$ and $n+1$ to the corresponding subalgebra of
$k[{\Delta_1}]$; see Lemma \ref{vero} below. Taking normalizations yields an isomorphism $\bar B\approx
k[{c\Delta_1}]$. But then $B=k[{c\Delta_1}]$ as well, because $\bar B$ and $B$
coincide in degree 1 (being retracts of algebras with this property).
\end{proof}

\begin{lemma}\label{vero}
Let $A$ and $B$ be $\ZZ$-graded rings. Suppose that $B$ is reduced. If the homogeneous homomorphisms $\phi:A_{(n)}\to B_{(n)}$ and $\psi:A_{(n+1)}\to B_{(n+1)}$ coincide on $A_{(n(n+1))}$, then they have a common extension to a homogeneous homomorphism $\chi:V\to B$, where $V$ is the subalgebra of $A$ generated by  $A_{(n)}$ and $A_{(n(n+1))}$. If, in addition, $A$ is reduced and $\phi$ and $\psi$ are injective, then $\chi$ is also injective.
\end{lemma}

\begin{proof}
One checks easily that one only needs to verify the following: if
$uv=u'v'$ for homogeneous elements $u,u'\in A_{(n)}$, $v,v'\in
A_{(n+1)}$, then $\phi(u)\psi(v)=\phi(u')\psi(v')$. As $B$ is reduced, it is
enough that $(\phi(u)\psi(v)-\phi(u')\psi(v'))^{n(n+1)}=0$. Since
$$
u^p(u')^{n(n+1)-p},\ v^p(v')^{n(n+1)-p}\in A_{(n(n+1))}, \qquad p\in [0,n(n+1)],
$$
this follows immediately from the hypothesis that $\phi$ and $\psi$ coincide on $A_{(n(n+1))}$.

If $A$ is reduced, then every non-zero homogeneous ideal in $A$ intersects $A_{(n(n+1))}$ non-trivially, and this implies the second assertion.
\end{proof}

\section{Tensor factors}\label{tensor}

In this section we examine a special case of conjecture A. Suppose
$$
A=k\oplus A_1\oplus A_2\oplus\cdots\ \ \text{and}\ \ A=k\oplus
B_1\oplus B_2\oplus\cdots.
$$
are (standard) graded. Then, as usual, their tensor product (over $k$) is the
(standard) graded algebra
$$
A\otimes_k B=k\oplus C_1\oplus C_2\oplus\cdots,\ \
C_h=\bigoplus_{i+j=h}(A_i\otimes_k B_j).
$$
Clearly, $A$ is a graded retract of $A\otimes B$ -- one just considers
the retraction
$$
A\otimes B\>{1_A\otimes\pi_B}>>A\otimes k=A,
$$
where $\pi_B:B\to k$ is the augmentation with $\Ker(\pi_B)$ the
irrelevant maximal ideal.

In the special case of polytopal algebras $A=k[P]$ and $B=k[Q]$ we
have the following easily verified formula
$$
k[P]\otimes k[Q]\approx k[{\join(P,Q)}].
$$
We recall from \cite{BG2} that a lattice polytope $R$ is a join of two lattice polytopes $P$ and
$Q$ if $R$ has two faces isomorphic to $P$ and $Q$ as lattice polytopes, the
affine hulls of these faces do not intersect each other in the affine hull of
$R$, and $\L_R=\L_P\cup\L_Q$. Clearly, in this situation $R$ is the convex
hull of these faces. (Note that $R$ is unique up to an isomorphism of lattice
polytopes.)

\begin{theorem}\label{3.1}
If $A$ and $B$ are two standard graded $k$-algebras such that $A\otimes
B\approx k[P]$ for some normal lattice polytope $P$, then both $A$ and
$B$ are polytopal.
\end{theorem}

In particular, we obtain the {\em splitting off property} of polytopal
extensions:
$$
\bigl(A\otimes k[Q]\approx k[P]\bigr)\Rightarrow(A\approx k[R]\ \text{for
some}\ R)
$$
if $P$ is a normal lattice polytope.

An essential step in the proof of Theorem \ref{3.1} is the special case when
one of the factors collapses into the cone over a point, that is a
polynomial ring in one variable.

\begin{proposition}\label{3.2}
Let $A$ be as in Theorem \ref{3.1} and $P$ be any lattice polytope. If
$A[X]\approx k[P]$ as graded algebras ($X$ a variable) then there is a
polytope $R$ for which $A\approx k[R]$.
\end{proposition}

\begin{remark}\label{3.3}
Let us call an ideal $I$ {\em binomial} if it is generated by monomials $\alpha
T$ and binomials $T-\beta T'$ where $\alpha,\beta\in k^*$ and $T$, $T'$ denote
terms. The arguments in the proof of Proposition \ref{3.2}, presented below,
yield the following more general result: for a not necessarily algebraically
closed field $k$ and an ideal $I\subset k[X_1,\ldots,X_n]$ the ideal
$$
Ik[X_1,\ldots,X_n,X]\subset k[X_1,\ldots,X_n,X]
$$
can be transformed into a binomial ideal by a linear automorphism of
$k[X_1,\ldots,X_n,\allowbreak X]$ if and only if this is possible for $I$ in
$k[X_1,\dots,X_n]$.
\end{remark}

Several times we will use the following theorem \cite[2.6]{ES} characterizing
binomial prime ideals in affine semigroup rings over algebraically closed
fields. (In \cite{ES} the theorem is given only for polynomial rings, but the
generalization is immediate.)

\begin{theorem}\label{EisStu}
Let $k$ be an algebraically closed field. A binomial ideal $I$ in an affine
semigroup ring $k[S]$ is prime if and only if the residue class ring $k[S]/I$
contains a (multiplicative) affine semigroup $S'$ such that $k[S]/I=k[S']$ and,
moreover, the natural epimorphism $k[S]\to k[S']$ maps the monomials in $k[S]$
to those in $k[S']$.
\end{theorem}

\begin{proof}[Proof of Proposition \ref{3.2}]
Since polytopal algebras are precisely the standard graded affine
semigroup rings coinciding with their normalizations in degree 1
\cite{BGT}, we only need to show that $A$ is a semigroup ring.

Theorem \ref{EisStu} reduces the proposition to the assertion on binomial ideals
$I$ discussed in Remark \ref{3.3} where the ideal $I$ to be considered is the
kernel of a surjection $k[X_1,\dots,X_n]\to A$.

We set $X_{n+1}=X$. By hypothesis there exists a linear transformation $\Psi$,
$$
\Psi(X_i)=t_{i1}X_1+\cdots+t_{in}X_n+t_{in+1}X_{n+1},\qquad i\in[1,n+1],\quad t_{ij}\in k,
$$
such that $\Psi(Ik[X_1,\dots,X_{n+1}])$ is binomial. For all $j\in[1,n]$ we
consider the surjection
$$
\pi_j: k[X_1,\dots,X_{n+1}]\to k[X_1,\dots,X_n],\quad \pi_j(X_i)=X_i,\ \pi_j(X_{n+1})= X_j.
$$

We define the linear endomorphisms $\epsilon_0$ and $\epsilon_j$, $j\in[1,n]$ of
$k[X_1,\dots,X_n]$ by
\begin{align*}
\epsilon_0(X_i)&=t_{i1}X_1+\cdots+t_{in}X_n,\ \ \ i\in[1,n]\\
\epsilon_j(X_i)&=t_{i1}X_1+\cdots+t_{in}X_n+ t_{in+1}X_j,\ \ \ i\in[1,n].
\end{align*}

There exists $j$ for which $\epsilon_j$ is an automorphism. In fact, let $T$
denote the $(n+1)\times(n+1)$ matrix $(t_{ij})$. If $\epsilon_0$ is not an
automorphism, then the $n\times n$ submatrix $T'$ of $T$ consisting of the
intersection of the first $n$ rows with the first $n$ columns of $T$ has rank
$n-1$ and the larger submatrix $T''$ consisting of the first $n$ rows has rank
$n$ (otherwise $\rank T<n+1)$.) There exists $j\in[1,n]$ such that the $j$th
column of $T'$ lies in the vector space $U$ spanned by the remaining $n-1$
columns of $T'$, but the $(n+1)$th column of $T''$ does not belong to $U$. If we
replace the $j$th column of $T'$ by its sum with $(n+1)$st column of $T''$, we
obtain an $n\times n$ matrix of rank $n$, and this matrix defines the
automorphism $\epsilon_j$.

There exists a unique epimorphism $\nu:k[X_1,\ldots,X_n,X_{n+1}]\to
k[X_1,\ldots,\allowbreak X_n]$ making the diagram
$$
\begin{CD}
k[X_1,\ldots,X_n,X_{n+1}]@>\Psi>>k[X_1,\ldots,X_n,X_{n+1}]\\
@V\nu VV@VV\pi_j V\\
k[X_1,\ldots,X_n]@>\epsilon_j>>k[X_1,\ldots,X_n],
\end{CD}
$$
commutative. Since $\nu(X_i)=X_i$ for $i\in[1,n]$, we have
$\nu(Ik[X_1,\ldots,X_n,\allowbreak X_{n+1}])=I$, and $\pi_j$ maps binomial
ideals to binomial ideals.
\end{proof}

Two more auxiliary facts:

\begin{lemma}\label{3.4}
Let $S$ be a positive affine semigroup. Then the graded ring
$$
\gr_{\mathfrak m}(k[S])=\gr_{\mathfrak m}(k[S]_{\mathfrak m})
$$
is a domain if and only if $S$ is homogeneous.
\end{lemma}

\begin{proof}
If $S$ is not homogeneous, then among the relations between the irreducible
elements of $S$ there occurs an equation
$$
s_1^{a_1}\cdots s_r^{a_r}=s_1^{b_1}\cdots s_r^{b_r},\quad
a_1+\cdots+a_r<b_1+\cdots+b_r,
$$
for some $a_i,b_i\geq0$. But then
$\overline{s_1}^{a_1}\cdots\overline{s_r}^{a_r}=0$ in $\gr_{\mathfrak
m}(k[S])$, where $\bar s_i$ denotes the residue class of $s_i$ in
${\mathfrak m}/{\mathfrak m}^2$. Since these are nonzero classes we get a
contradiction.

The other implication is clear.
\end{proof}

\begin{lemma}\label{3.5}
Let $S$ be a normal positive affine semigroup. For an edge $E$ of the polyhedral
cone $C(S)$ spanned by $S$ in $\RR\otimes\gp(S)$ let $x_E$ be the generator of
the subsemigroup $E\cap S\subset S$ (isomorphic to $\ZZ_+$) and let $x_E^{-1}S$
denote the sub-semigroup of $\gp(S)$ generated by
$S\cup\{x_E^{-1},x_E^{-2},\ldots\}$. Then
$$
x_E^{-1}S\approx\ZZ\times S'
$$
for some normal positive affine semigroup $S'$.
\end{lemma}
For example, this is proved in \cite[\S2]{Gu}.

\begin{proof}[Proof of Theorem \ref{3.1}.]
As $A$ and $B$ are retracts of $k[P]$
they have to be normal domains.

We denote by $X$ and $Y$ the affine normal toric varieties corresponding
to $A$ and $B$. The two irrelevant maximal ideals will be denoted by
${\mathfrak m}_A$ and ${\mathfrak m}_B$ respectively.

Let $y\in Y$ be any smooth point. Then for the local ring of the
point $({\mathfrak m}_A,y)\in X\times Y$ we have
$$
\gr({\Oo}_{({\mathfrak m}_A,y)})\approx A[Y_1,\ldots,Y_{\dim(B)}].
$$
as graded algebras (the $Y_i$ are indeterminates over $A$).

On the other hand, identifying $X\times Y$ with the variety
$\maxSpec(k[P])$ we see that there is a maximal ideal ${\mathfrak n}\subset
k[P]$ such that
$$
\gr(\Oo_{(\mathfrak m_A,y)})\approx\gr_{\mathfrak n}(k[P]_{\mathfrak n}).
$$
By an iterated use of Lemma \ref{3.5} we get
$$
k[P]_{\mathfrak n}\approx(k[S][T_1,T_1^{-1},\ldots,T_s,T_s^{-1}])_{ \mathfrak r}
$$
for some positive normal affine semigroup $S$, an integer $s\geq0$, and a
maximal ideal $ \mathfrak r\subset k[S][T_1,T_1^{-1},\ldots,T_s,T_s^{-1}]$,
containing $S\setminus\{1\}$. Since $k$ is algebraically closed there are
elements $a_1,\ldots,a_s\in k^*$ such that $ \mathfrak r$ is generated by the
maximal monomial ideal $ \mathfrak r_0\subset k[S]$ and
$\{T_1-a_1,\ldots,T_s-a_s\}$. Therefore one has a natural graded isomorphism
$$
\gr_{{\mathfrak n}}(k[\bar S_P]_{{\mathfrak n}})\approx(\gr_{ \mathfrak r_0}(k[S]_{ \mathfrak r_0}))
[\tau_1,\ldots,\tau_s]
$$
for $\tau_i=T_i-a_i$.

Consequently,
$$
A[Y_1,\ldots,Y_{\dim(B)}]\approx(\gr_{ \mathfrak r_0}(k[S]_{ \mathfrak r_0}))[\tau_1,
\ldots,\tau_s]
$$
as graded algebras. By Lemma \ref{3.4} $S$ is a homogeneous affine semigroup.
Then, clearly, $S\approx S_Q$ for some normal polytope $Q$ and we get the
graded isomorphism
$$
A[Y_1,\ldots,Y_{\dim(B)}]\approx k[Q][\tau_1,
\ldots,\tau_s]\approx k[{\join(Q,\Delta_{s-1}]}]
$$
($\Delta_{s-1}$ denotes the $(s-1)$-unit simplex). By an
iterated use of Proposition \ref{3.2} we can split off the polynomial extension
on the left. That is, $A$ is polytopal. By symmetry $B$ is also polytopal.
\end{proof}

\begin{remark}
As a consequence of \ref{3.2} one can show that any retract of a polytopal
algebra defined by a single equation is polytopal. In fact, suppose that
$k[P]=k[X_1,\dots,X_n]/(b)$ where $b$ is a binomial of degree $d>0$ and denote
the natural epimorphism from $k[X_1,\dots,X_n]$ onto $k[P]$ by $\pi$. Let $h$ be
a retraction of $k[P]$, $y_1,\dots,y_m$ a $k$-basis of $\Im(h)_1$ and
$z_1,\dots,z_p$ a $k$-basis of $\Ker(h)_1$. Let $g:k[Y_1,\dots,Y_m]\to \Im(h)$
be given by $g(Y_i)=y_i$, and let $\gamma$ be the (unique) epimorphism
$k[X_1,\dots,X_n]\to k[Y_1,\dots,Y_m]$ with $g\circ \gamma=h\circ \pi$. Then
$\Ker(g)$ is generated by $\pi(b)$. If $\pi(b)=0$, then $\Im(h)$ is a polynomial
ring. Otherwise $b'=\pi(b)$ is of degree $d$. Now we extend $g$ to the
epimorphism $g':(k[Y_1,\dots,Y_m]/(b'))[Z_1,\dots,Z_p]\to k[P]$ with
$g'(Z_i)=z_i$. Comparing Hilbert functions we see that $g'$ is an isomorphism,
and Proposition \ref{3.2} shows that $\Im(h)\approx k[Y_1,\dots,Y_m]/(b')$ is
polytopal.
\end{remark}

\section{Polytopal linear groups}\label{aut}

Here we survey the relevant part of \cite{BG1} that will be used in the
sequel. Throughout this section $k$ is a general (not necessarily
algebraically closed) field.

For a lattice polytope $P$ the group $\Gamma_k(P)=\ggr.aut(k[P])$ is a
linear $k$-group in a natural way. It coincides with $\GL_n(k)$ in the
case $P$ is the unit $(n-1)$-simplex $\Delta_{n-1}$. The groups
$\Gamma_k(P)$ are called {\em polytopal linear groups} in \cite{BG1}.

An element $v\in \ZZ^d$, $v\neq 0$, is a {\em column vector} (for $P$) if
there is a facet $F\subset P$ such that $x+v\in P$ for every
lattice point $x\in P\setminus F$. The facet $F$ is called {\em the base
facet} of $v$.
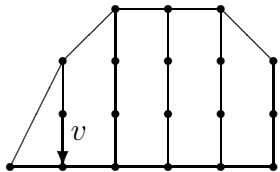
\begin{figure}[htb]
\begin{center}
\begin{picture}(5,3)(0,0)
\multiput(0,0)(1,0){6}{\vertex}
\multiput(1,1)(1,0){5}{\vertex}
\multiput(1,2)(1,0){5}{\vertex}
\multiput(2,3)(1,0){3}{\vertex}
\put (0,0){\line(1,0){5}}
\put (0,0){\line(1,2){1}}
\put (1,2){\line(1,1){1}}
\put (2,3){\line(1,0){2}}
\put (4,3){\line(1,-1){1}}
\put (5,2){\line(0,-1){2}}
\put (1,0){\line(0,1){2}}
\put (2,0){\line(0,1){3}}
\put (3,0){\line(0,1){3}}
\put (4,0){\line(0,1){3}}
\put (1.15,0.5){$v$}
\thicklines
\put (1,1){\vector(0,-1){1}}
\end{picture}
\end{center}
\caption{A column structure}\label{Fig2}
\end{figure}
The set of column vectors of $P$ is denoted by $\Col(P)$. A pair $(P,v)$,
$v\in\Col(P)$, is called a {\em column structure}. Let $(P,v)$ be a column
structure and $P_v\subset P$ be the base facet for $v\in\Col(P)$. Then for each
element $x\in S_P$ we set $\het_v(x)=m$ where $m$ is the largest non-negative
integer for which $x+mv\in S_P$. Thus $\het_v(x)$ is the `height' of $x$ above
the facet of the cone $C(S_P)$ corresponding to $P_v$ in direction $-v$. It is
an easy observation that $x+\het_v(x)\cdot v\in S_{P_v}\subset S_P$ for any
$x\in S_P$.

Let $(P,v)$ be a column structure and $\lambda\in k$. We identify the
vector $v$, representing the difference of two lattice points in $P$,
with the degree $0$ element $(v,0)\in\gp(S_P)\subset k[\gp(S_P)]$. Then
the assignment
$$
x\mapsto (1+\lambda v)^{\het_v(x)}x.
$$
gives rise to a graded $k$-algebra automorphism $e_v^\lambda$ of $k[P]$. Observe
that $e_v^\lambda$ becomes an elementary matrix in the special case when
$P=\Delta_{n-1}$, after the identifications
$k[{\Delta_{n-1}}]=k[X_1,\ldots,X_n]$ and $\Gamma_k(P)=\GL_n(k)$. Accordingly
$e_v^\lambda$ is called an {\em elementary automorphism}.

In the following $\AA_k^s$ denotes the additive group of the $s$-dimensional affine
space.

\begin{proposition}\label{6.1}
Let $v_1,\ldots,v_s$ be pairwise different column vectors for $P$ with
the same base facet $F=P_{v_i}$, $i=1,\dots,s$.
\begin{itemize}
\item[(a)]
The mapping
$$
\phi:\AA_k^s\to\Gamma_k(P),\qquad
(\lambda_1,\ldots,\lambda_s)\mapsto
e_{v_1}^{\lambda_1}\circ\cdots\circ e_{v_s}^{\lambda_s},
$$
is an embedding of algebraic groups. In particular,
$e_{v_i}^{\lambda_i}$ and $e_{v_j}^{\lambda_j}$ commute for all
$i,j\in\{1,\dots,s\}$ and $(e_{v_1}^{\lambda_1}\circ\cdots\circ
e_{v_s}^{\lambda_s})^{-1} =e_{v_1}^{-\lambda_1}\circ\cdots\circ
e_{v_s}^{-\lambda_s}.$ \item[(b)]
For $x\in\L_P$ with $\het_{v_1}(x)=1$ one has
$$
e_{v_1}^{\lambda_1}\circ\cdots\circ e_{v_s}^{\lambda_s}(x)=(1+
\lambda_1v_1+\cdots+\lambda_sv_s)x.
$$
\end{itemize}
\end{proposition}

The image of the embedding $\phi$ given by Lemma \ref{6.1} is denoted by
$\AA(F)$. Of course, $\AA(F)$ may consist only of the identity map of
$k[P]$, namely if there is no column vector with base facet $F$.

Put $n=\dim(P)+1$. The $n$-torus $\TT_n=(k^*)^n$ acts naturally on
$k[P]$ by restriction of its action on $k[\gp(S_P)]$ that is given by
$$
(\xi_1,\ldots,\xi_{n+1})(e_i)=\xi_ie_i,\quad
i\in[1,n].
$$
Here $e_i$ is the $i$-th element of a fixed basis of
$\gp(S_P)=\ZZ^n$.  This gives rise to an algebraic embedding
$\TT_n\subset\Gamma_k(P)$, whose image we denote by $\TT_k(P)$.
It consists precisely of those automorphisms of $k[P]$ which
multiply each monomial by a scalar from $k^*$.

The (finite) automorphism group $\Sigma(P)$ of the semigroup $S_P$ is also a
subgroup of $\Gamma_k(P)$.

\begin{theorem}\label{6.2}
Let $P$ be a convex lattice $n$-polytope and $k$ a field.
Every element $\gamma\in\Gamma_k(P)$ has a (not uniquely determined)
presentation
$$
\gamma=\alpha_1\circ\alpha_2\circ\cdots\circ\alpha_r\circ\tau\circ\sigma,
$$
where $\sigma\in\Sigma(P)$, $\tau\in\TT_k(P)$, and
$\alpha_i\in\AA(F_{i})$ such that the facets $F_i$ are pairwise
different and $\#L_{F_i}\le \#L_{F_{i+1}}$, $i\in[1,r-1]$.
\end{theorem}

For any positive affine semigroup $S$ and any grading
$$
k[S]=k\oplus A_1\oplus A_2\oplus\cdots,
$$
making all non-unit monomials homogeneous of positive degree, the
automorphisms of $k[S]$ that multiply monomials by scalars from $k^*$
constitute a closed $\rank(S)$-torus
$$
\TT_k(S)\subset\ggr.aut_k(k[S]),\ \ \ \TT_k(S)\approx(k^*)^{\rank(S)}.
$$
This torus is naturally identified with the set of $k$-rational points of the open subscheme
$$
\Spec(k[\gp(S)])\subset\Spec(k[S]).
$$
and is called the {\em embedded torus}. In case $k$ is
infinite, the embedded torus is maximal (by literally the same arguments as at
the end of \S3 in \cite{BG1}).

\section{The structure of retractions}\label{struct}

Now we first consider Conjecture B in detail and then observe that it
does not admit a direct extension to codimension $\ge 2$.

Let $P\subset\RR^n$ be a lattice polytope of dimension $n$ and $F\subset
P$ a face. Then there is a uniquely determined retraction
$$
\pi_F:k[P]\to k[F],\ \ \pi_F(x)=0\quad \text{for}\ x\in\L_P
\setminus F.
$$
Retractions of this type will be called {\em face retractions} and {\em facet retractions} if $F$ is a facet or, equivalently, $\codim(\pi_F)=1$.

Now suppose there are an affine subspace $H\subset\RR^n$
and a vector subspace $W\subset\RR^n$ with $\dim W+\dim H=n$, such that
$$
\L_P\subset\bigcup_{x\in\L_P\cap H}(x+W).
$$
(Observe that $\dim(H\cap P)=\dim H$.) The triple $(P,H,W)$ is called a {\em
lattice fibration of codimension $c=\dim W$}, whose {\em base polytope} is
$P\cap H$; its {\em fibers} are the maximal lattice subpolytopes of $(x+W)\cap
P$, $x\in\L_P\cap H$ (the fibers may have smaller dimension than $W$). $P$
itself serves as a {\em total polytope} of the fibration. If $W=\RR w$ is a
line, then we call the fibration {\em segmental} and write $(P,H,w)$ for it.
\begin{figure}[htb]
\begin{center}
\begin{picture}(5,5)(0,0)
\multiput(3,0)(1,0){3}{\vertex}
\multiput(1,1)(1,0){5}{\vertex}
\multiput(1,2)(1,0){5}{\vertex}
\multiput(2,3)(1,0){3}{\vertex}
\put (3,0){\line(1,0){2}}
\put (3,0){\line(-2,1){2}}
\put (1,2){\line(1,1){1}}
\put (2,3){\line(1,0){2}}
\put (4,3){\line(1,-1){1}}
\put (5,2){\line(0,-1){2}}
\put (1,1){\line(0,1){1}}
\put (2,0.5){\line(0,1){2.5}}
\put (3,0){\line(0,1){3}}
\put (4,0){\line(0,1){3}}
\put (0,1){\line(1,0){6}}
\put (5.5,1.2){$H$}
\put (2.2,1.2){$w$}
\thicklines
\put (2,1){\vector(0,1){1}}
\end{picture}
\end{center}
\caption{A lattice segmental fibration}\label{Fig1}
\end{figure}
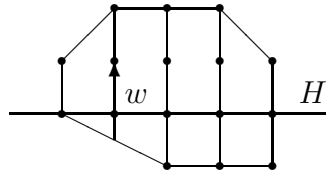
Note that the column structures introduced in Section \ref{aut} give rise to
lattice segmental fibrations in a natural way.

For a lattice fibration $(P,H,W)$ let $L\subset\ZZ^n$ denote the subgroup
spanned by $\L_P$, and let $H_0$ be the translate of $H$ through the
origin. Then one has the direct sum decomposition
$$
L=(L\cap W)\oplus(L\cap H_0).
$$
Equivalently,
$$
\gp(S_P)=L\oplus\ZZ= \bigl(\gp(S_P)\cap W_1\bigr)\oplus
\gp(S_{P\cap H_1}).
$$
where $W_1$ is the image of $W$ under the embedding $\RR^n\to\RR^{n+1}$,
$w\mapsto (w,0)$, and $H_1$ is the vector subspace of $\RR^{n+1}$ generated by
all the vectors $(h,1)$, $h\in H$.

For a fibration $(P,H,W)$ one has the naturally associated
retraction:
$$
\rho_{(P,H,W)}:k[P]\to k[{P\cap H}];
$$
it maps $\L_P$ to
$\L_{P\cap H}$ so that fibers are contracted to their
intersection points with the base polytope $P\cap H$.
\medskip

Clearly, if $f:k[P]\to k[P]$ is a retraction, then for any graded automorphism
$\alpha$ of $k[P]$ the composite map $f^{\alpha}= \alpha\circ f\circ\alpha^{-1}$
is again a retraction and $\Im(f^{\alpha})=\alpha(\Im(f))$ and
$\Ker(f^{\alpha})=\alpha(\Ker(f))$. Now the exact formulation of Conjecture B
is as follows. \medskip

\noindent{\bf Conjecture B.} {\em For a codimension 1 retraction $f:k[P]\to k[P]$
there is $\alpha\in\Gamma_k(P)$ such that $f^{\alpha}=\iota\circ g$ for
a retraction $g$ of type either $\pi_F$ or $\rho_{(P,H,w)}$ and
$\iota:\Im(g)\to k[P]$ a graded $k$-algebra embedding.}
\medskip

In other words this conjecture claims that any codimension 1 retraction
can be `modified' by an automorphism so that the corrected retraction
factors through a retraction preserving the monomial structure.

The explicit description of the embeddings $\iota$
is not addressed in this paper.

A necessary condition for Conjecture B is that any codimension $1$ retraction $f$ can be
modified by a graded automorphism $\alpha$ so that
$f^{\alpha}$ has either a homogeneous binomial of degree $1$ or a monomial of
degree 1 in its kernel. A weaker condition is that $\Ker(f^\alpha)$ contains a
homogeneous binomial of degree $\ge 1$ (evidently this holds if there is a monomial
in $\Ker(f^\alpha)$).

We remark that even an example of just an endomorphism in
$\Pol_k$, such that $\Ker(f^{\alpha})$ contains no (homogeneous)
binomial for any $\alpha$, is not readily found. However, such
exists, even in the class of codimension 2 retractions.

The examples below are constructed from joins of polytopes. The following lemma
enables us to describe $\Gamma_k(R)$ for $R=\join(P,Q)$ under a mild assumption
on $P$ and $Q$. We identify them with the corresponding faces of $R$.

\begin{lemma}\label{4.1}
Let $P$ and $Q$ be lattice polytopes, both having interior lattice
points. Then $\Col(\join(P,Q))=\Col(P)\cup\Col(Q)$.
\end{lemma}

\begin{proof}
That each of the column vectors of the polytopes serves as a column
vector for $\join(P,Q)$ is clear.

Now let $v\in\Col(\join(P,Q))$. If $v$ is parallel to either $P$ or $Q$ then
either $v\in\Col(P)$ or $v\in\Col(Q)$ since $\L_{\join(P,Q)}=\L_P\cup\L_Q$. So
without loss of generality we can assume that $v$ is parallel neither to $P$ nor
to $Q$. Since $P$ and $Q$ span $\join(P,Q)$ they cannot be contained
simultaneously in the base facet of $v$. But then either $p+v\in\join(P,Q)$ or
$q+v\in\join(P,Q)$ for suitable vertices $p\in P$ and $q\in Q$. We get a
contradiction because one of the points $r+v$ or $s+v$ is outside $\join(P,Q)$ for
interior lattice points $r\in P$, $s\in Q$.
\end{proof}

\begin{example}\label{4.2}
Let $Q$ be the lattice triangle spanned by $(0,-1)$, $(-1,0)$, and $(1,1)$. Then
$Q$ contains only one more lattice point, namely $(0,0)$. Identifying $U$ with
$(0,0)$, $V$ with $(0,-1)$, and $W$ with $(-1,0)$ we see that the polynomial
ring $k[U,V,W]$ can be embedded into $k[Q]$ such that the indeterminates
correspond to lattice points. Moreover, $k[\gp(S_Q)]$ is then just the Laurent
polynomial ring $k[\ZZ^3]=k[U^{\pm 1},V^{\pm 1},W^{\pm 1}]$.

Let $h':k[X,Y]\to k[U,V,W]$ be defined by $h'(X)=U+V$, $h'(Y)=U+W$. Then $h'$
induces a retraction $h$ of $k[U,V,W,X,Y]$, namely the retraction mapping $X$
and $Y$ to $h'(X)$ and $h'(Y)$ respectively and leaving $U,V,W$ invariant. This
retraction extends in a natural way to retraction of $k[U^{\pm 1},V^{\pm
1},W^{\pm 1},X,Y]$, and can then be restricted to
$$
k[Q]\tensor k[\Delta_1]\subset k[U^{\pm 1},V^{\pm 1},W^{\pm 1},X,Y]
$$
where we identify $k[X,Y]$ with the polytopal ring $k[\Delta_1]$ of the unit
segment. It can further be restricted to $k[\join(2Q,2\Delta_1)]$ which is
embedded into $k[Q]\tensor k[\Delta_1]$ as the tensor product of the second
Veronese subalgebras of the normal algebras $k[Q]$ and $k[\Delta_1]$.

We claim that the just constructed retraction $h$ of $k[P]$,
$P=\join(2Q,2\Delta_1)$, $\dim P=4$, has no conjugate $h^\alpha$ by an automorphism
$\alpha\in\Gamma_k(P)$ such that the kernel of $h^\alpha$ contains a binomial.

The polytope $Q$ has no column structures, a property inherited by $2Q$.
Moreover, both $2Q$ and $2\Delta_1$ have interior points. Therefore the only
column structures on $P$ are those it gets from $2\Delta_1$ (see Lemma
\ref{4.1}). Then every element $\alpha\in\Gamma_k(P)$ is of the form
$\tau\circ\beta$, where $\tau$ is a toric automorphism and
$\beta=1\tensor\beta'$ for some $\beta'\in\Gamma_k(2\Delta_1)$. Since $\tau$
does not affect the monomial structure, we can assume $\tau=1$. Furthermore the
graded automorphisms of $k[2\Delta_1]$ are all restrictions of automorphisms of
$k[\Delta_1]=k[X,Y]$ so that we have to take into account all automorphisms of
$k[P]$ induced by a substitution
$$
X\mapsto a_{11}X+a_{12}Y,\quad Y\mapsto a_{21}X+a_{22}Y,\quad U\mapsto U,\quad
V\mapsto V,\quad W\mapsto W
$$
with $\det(a_{ij})\neq 0$. Then $h^\alpha$ is induced by the substitution
$$
a_{11}X+a_{12}Y\mapsto U+V,\qquad a_{21}X+a_{22}Y\mapsto U+W,
$$
leaving $U,V,W$ invariant. Also $h^\alpha$ extends to a retraction of $k[U^{\pm
1},V^{\pm 1},\allowbreak W^{\pm 1},\allowbreak X,Y]$ and then restricts to $k[U,V,W,X,Y]$. This shows that
the kernel of the extension cannot contain a monomial; otherwise it would
contain a monomial in $X$ and $Y$, but $h^\alpha$ is injective on $k[X,Y]$. If the kernel contains a binomial $b$, we can
assume that $b\in k[U,V,W,X,Y]$. In other words, we can find a binomial in the
ideal $\pp$ of $k[U,V,W,X,Y]$ generated by
$$
a_{11}X+a_{12}Y-(U+V),\qquad a_{21}X+a_{22}Y-(U+W).
$$
Since the prime ideal $\pp$ contains no monomials, we can assume that the two
terms of $b$ are coprime. But then $b$ reduces to a monomial modulo one of the
variables, and since $\pp$ reduces to an ideal generated by linear forms, it reduces to a prime ideal. The reduction of $\pp$ modulo any of the variables cannot contain another variable.
\end{example}

In view of what has been said above it is natural to introduce the following
classes of {\em tame retractions} (for arbitrary codimension) and {\em tame
morphisms} in $\Pol_k$: a retraction $f:k[P]\to k[P]$ is called
tame if there are a sequence of lattice polytopes
$P=P_m,P_{m-1}\ldots,P_1,P_0=Q$, automorphisms $\alpha_i\in\Gamma_k(P_i)$ and
monomial structure preserving retractions $\pi_i:k[{P_i}]\to k[{P_{i-1}}]$
such that
$$
f=\iota \circ\pi_1\circ\alpha_1\circ\cdots\circ\pi_{m-1}\circ\alpha_{m-1}
\circ\pi_m\circ\alpha_m
$$
for a $k$-algebra embedding $\iota:k[Q]\to k[P]$.

The definition of a tame morphism is literally the same with the word
`retraction' changed to `morphism'. Notice that the image of a tame
retraction is again a polytopal ring, whereas that of a tame morphism
is a semigroup ring of some homogeneous semigroup.

In this terminology Conjecture B merely says that all codimension 1
retractions are tame, while Example \ref{4.2} shows that there are polytopal
codimension 2 retractions, which are not even tame morphisms. The following example gives a wild retraction that as a morphism is tame.

\begin{example}\label{WildTame}
Let $P$ be the join of $7\Delta_1$ and $2\Delta_1$, and thus $\dim P=3$ (see
Figure \ref{FigTame}).
\begin{figure}[htb]
\begin{center}
\begin{picture}(7.5,3.5)(0,-1.5)
\multiput(0,2)(1,0){8}{\vertex}
\multiput(3.5,-1.5)(0,0.875){3}{\vertex}
\put (0,2){\line(1,0){7}}
\put (3.5,-1.5){\line(0,1){1.75}}
\put (0,2){\line(2,-1){3.5}}
\put (0,2){\line(1,-1){3.5}}
\put (7,2){\line(-2,-1){3.5}}
\put (7,2){\line(-1,-1){3.5}}
\put (-0.7,1.8){$A_1$}
\put (7.3,1.8){$A_8$}
\put (3.8,-1.7){$B_1$}
\put (3.3,0.6){$B_3$}
\end{picture}
\quad
\begin{picture}(7.5,3.5)(0,-0.5)
\multiput(0,0)(1,0){8}{\vertex}
\multiput(0,1)(1,0){5}{\vertex}
\put (0,0){\line(1,0){7}}
\put (0,1){\line(1,0){4}}
\put (0,0){\line(0,1){1}}
\put (4,1){\line(3,-1){3}}
\put (-0.7,-0.2){$A_1'$}
\put (7.3,-0.2){$A_8'$}
\put (-0.2,1.2){$B_1'$}
\put (1.8,1.2){$B_2'$}
\put (3.8,1.2){$B_3'$}
\thicklines
\put (1,1){\vector(1,-1){1}}
\put (1,1){\vector(2,-1){2}}
\end{picture}
\end{center}
\caption{The polytopes $P$ and $P'$}\label{FigTame}
\end{figure}
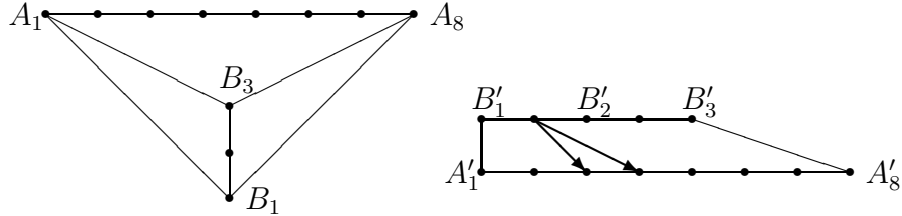
It is straightforward to check that the assignment
$$
B_1\mapsto A_2+A_3,\quad
B_2\mapsto A_4+A_5,\quad
B_3\mapsto A_6+A_7, \quad A_i\mapsto A_i,\quad i\in[1,8]
$$
defines a retraction $h$ of $k[P]$, and it is easily checked that $h=\iota\circ
g\circ f$, where $\iota$ is the natural embedding of $k[7\Delta_1]$ into $k[P]$,
$f$ is the morphism $k[P]\to k[P']$ given by the assignment $A_i\mapsto A_i'$,
$B_i\mapsto B_i'$, and $g:k[P']\to k[7\Delta_1]$ is a tame surjection. In fact,
one can choose $g=\pi_{7\Delta_1}\circ\alpha$, where $\alpha\in\AA(7\Delta_1)$
is an automorphism of $k[P']$. (The column vectors used for $\alpha$ are
indicated in the figure.)

The kernel of $h$ is generated by $B_1-(A_2+A_3)$, $B_2-(A_4+A_5)$,
$B_3-(A_6+A_7)$. Assume there is an automorphism $\gamma\in\Gamma_k(P)$ such
that the image of this ideal contains a binomial or monomial of degree 1. Then
$$
b=\gamma(c_1(B_1-(A_2+A_3))+c_2(B_2-(A_4+A_5))+c_3(B_3-(A_6+A_7))
$$
is a binomial or monomial for suitable $c_1,c_2,c_3$, $c_i\neq 0$ for at least
one $i$. In view of the definition of $P$ as the join of two lattice segments
with interior points, we can write $\gamma$ as the composite of two
automorphisms induced by automorphisms $\gamma_1$ of $k[7\Delta_1]$ and
$\gamma_2$ of $k[2\Delta_1]$, respectively. In particular, if a polynomial
contains a term `on the edge $7\Delta_1$', then so does its image under
$\gamma$, and the same holds true for the edge $2\Delta_1$. Thus $b$ is not a monomial.
But neither can it be a binomial because then the polynomial
$$
\gamma_1(c_1(A_2+A_3)+c_2(A_4+A_5)+c_3(A_6+A_7))
$$
would be a term $a$ with
$$
\gamma_1^{-1}(a)=c_1(A_2+A_3)+c_2(A_4+A_5)+c_3(A_6+A_7).
$$
The right hand side is an `interior polynomial' of $k[7\Delta_1]$. On the other hand, an automorphism mapping a term to an `interior polynomial' must preserve the monomial structure by \cite[Lemma 4.1]{BG1}.
\end{example}

\section{Segmentonomial ideals}\label{seg}

For an affine semigroup $S$ an element $f\in k[S]$ will be called {\em
segmentonomial} if the Newton polytope $\N(f)\subset\RR\otimes\gp(S)$
has dimension $\leq 1$. (Clearly, monomials as well as binomials are
segmentonomials.) An ideal $I\subset k[S]$ is called segmentonomial if it
is generated by a system of segmentonomials.

It is proved in \cite{ES} that every minimal prime ideal over a binomial ideal
of $k[X_1,\ldots,X_n]$ (a polynomial ring) is again binomial. In this
section we derive the same result for
segmentonomial ideals in arbitrary affine semigroup rings.

\begin{theorem}\label{5.1}
Let $S$ be an affine semigroup and $I\subset k[S]$ be a segmentonomial
ideal.
\begin{itemize}
\item[(a)]
A minimal prime overideal $I\subset{\pp}\subset k[S]$ is binomial, and
$k[S]/I$ is again an affine semigroup ring.
\item[(b)]
Suppose that $\het(I)=1$, $f\in I$, $\dim(\N(f))=1$, and $\pp$ is as
above. Then for every system of pairwise distinct lattice points
$x_1,\ldots,x_m\in\L_P$, such that none of the pairs $(x_i,x_j)$,
$i\neq j$, spans a line in $\RR\otimes\gp(S) $parallel to $\N(f)$, the
residue classes $\bar x_1,\ldots,\bar x_m$ constitute a $k$-linearly
independent subset of $k[S]/\pp$.
\end{itemize}
\end{theorem}

\begin{proof}
We prove claim (a) by induction on $r=\rank(S)$. Claim (b) will follow
automatically from the description of $\pp$ derived below.

For $r=0$ there is nothing to show. Assume the theorem is proved for semigroups
of $\rank<r$ and choose a segmentonomial $f\in I$. Then $\pp$ contains a minimal
prime ${\pp}_0$ over the principal ideal $(f)$. Assume that ${\pp}_0$ is a
binomial ideal. By Theorem \ref{EisStu}, $k[S]/{\pp}_0\approx k[S_1]$ for some affine
semigroup $S_1$ and such that monomials in $k[S]$ go to monomials in $k[S_1]$.
But then segmentonomials in $k[S]$ are likewise mapped to segmentonomials in
$k[S_1]$. This holds true because affinely independent terms lift to affinely
independent terms. By induction hypothesis the image of $\pp$ in $k[S_1]$ is
binomial. Since binomials can be lifted to binomials in $k[S]$, we conclude that
$\pp$ is binomial.

The general situation thus
reduces to the case in which $I=(f)$ for some segmentonomial $f\in k[S]$ and $\het(I)=\het({\pp})=1$.

If $\pp$ contains
a monomial, then $\pp$ is a height 1
monomial prime ideal, and we are done.

Otherwise $S\cap{\pp}=\emptyset$. Consider the localization
${\pp}\,k[\gp(S)]$. It is a height 1 prime
ideal in the Laurent polynomial ring $k[\gp(S)]$. Therefore,
${\pp}\,k[\gp(S)]=f_0k[\gp(S)]$ for some prime element $f_0\in
k[\gp(S)]$.

Also $f_0$ is segmentonomial. In fact, we have
$f=f_0f_1$ for some $f_1\in k[\gp(S)]$ implying the equality
$
\N(f)=\N(f_0)+\N(f_1)
$
for the corresponding Newton polytopes. Since $\dim
N(f_0)=0$ is excluded, $\dim(\N(f_0))=1$. Multiplying $f_0$ by a suitable term
from $\gp(S)$ we can achieve that the origin $0\in\RR\otimes\gp(S)$ is one of
the end-points of $\N(f_0)$.

Let $\ell\subset\RR\otimes\gp(S)$ denote a rational line containing
$\N(f_0)$. In a suitable basis of the free abelian group $\gp(S)$ the
line $\ell$ becomes a coordinate direction. Therefore, we can assume that
$$
k[\gp(S)]=k[X_1,X_1^{-1},\ldots,X_n,X_n^{-1}]
$$
and that $f_0$ is a monic polynomial in $X_1$.
Since $k$ is algebraically closed, if follows  that
$
f_0=X_1-a
$
for some $a\in k$. Since $\pp$ does not contain a monomial, one has
$$
{\pp}=\pp\, k[\gp(S)]\cap k[S]=(X_1-a)k[X_1,X_1^{-1},\ldots,X_n,X_n^{-1}]\cap k[S].
$$
Thus $\pp$ is the kernel of the composite homomorphism
$$
k[S]\hookrightarrow k[X_1,X_1^{-1},\ldots,X_n,X_n^{-1}]
\>{X_1\mapsto a}>> k[X_2,X_2^{-1},\ldots,X_n,X_n^{-1}]
$$
This is a homomorphism mapping the elements of $S$ to Laurent monomials in $X_2,X_2^{-1},\ldots,X_n,X_n^{-1}$, and therefore $\pp$ is generated by binomials.
\end{proof}

\section{Based retractions}\label{base}

Throughout this section we suppose that $h:k[P]\to k[P]$ is a
retraction and that $A=\Im(h)$. We also assume
$P\subset\RR^n$, $\dim(P)=n$, $\gp(S_P)=\ZZ^{n+1}$ (and that $k$ is algebraically closed.)

\begin{lemma} \label{7.1}
The following conditions are equivalent:
\begin{itemize}
\item[(a)]
there is a subset $X\subset\L_P$ such that the restriction $h:k[S_X]\to A$ is an isomorphism, where $k[S_X]\subset
k[P]$ is the subalgebra generated by the semigroup $S_X=\langle X\rangle\subset S_P$,
\item[(b)]
there is a $(\dim(A)-1)$-dimensional cross section $Q$ of $P$
by a linear subspace $H$ such that $Q$ is a lattice polytope (i.~e. the
vertices of $Q$ are lattice points) and
$$
h|_{k[Q]}:k[Q]\to A
$$
is an isomorphism. In particular, $A$ is a polytopal algebra.
\end{itemize}
\end{lemma}
\begin{proof}
We only need to derive (b) from (a). Let $H$ be an affine hull of $X$ in
$\RR^n$. We have to show that $Q=H\cap P$ is a lattice polytope with $L_Q=X$.
Consider the subsemigroup
$$
S'_Q=\{x\in S|x\neq0\ \text{and}\ \RR_+x\cap P\subset H\}\cup\{0\}.
$$
Then $h(k[S'_Q])=A$ as well. On the other hand
$$
\dim k[S'_Q]=\dim H+1=\dim k[S_X]=\dim A.
$$
Thus the restriction $h:k[S'_Q]\to A$ is also an isomorphism. It follows that
$X=L_Q$, and every element in $S'_Q$ is a product of elements of $X$.
Furthermore $Q=\conv(X)$ since any rational point of the complement $Q\setminus
\conv(X)$ gives rise to elements in $S'_Q\setminus S_X$.
\end{proof}

A subpolytope $Q\subset P$ as in Lemma \ref{7.1}(b) (if it exists) will
be called a {\em base} of $h$ and $h$ is a {\em based
retraction}. Notice that a base is not necessarily uniquely determined.

\begin{theorem}\label{7.2a}
Suppose a retraction $h:k[P]\to k[P]$ has a base $Q$ that intersects the interior of $P$. Then
$h^{\tau}=\iota\circ\rho_{(P,H,W)}$ for some toric automorphism
$\tau\in\TT_k(P)$, a lattice fibration $(P,H,W)$ and a
$k$-algebra embedding $\iota:k[{H\cap P}]\to k[P]$.
In particular, $h$ is tame.
\end{theorem}

\begin{proof}
It is not hard to check that there is no restriction in assuming that $k[Q] = \Im(h)$.

Note that $\Ker(h)\cap S_P=\emptyset$. In fact, if a monomial is mapped to $0$
by $h$, then $\Ker(h)$ contains a monomial prime ideal $\pp$ of height $1$.
Since $\pp$ in turn contains all monomials in the interior of $S_P$, it must
also contain monomials from $S_Q$, which is impossible. Thus $h$ can be extended
to the normalization $k[\bar S_P]$; on $k[\bar S_Q]\subset k[\bar S_P]$ the
extension is the identity.

Set $L=\ZZ^{n+1}$, and let $U$ be the intersection of the $\QQ$-vector subspace
of $\QQ^{n+1}$ generated by $S_Q$ with $L$. Choose a basis $v_1,...,v_m$ of a
complement of $U$ in $L$. Since $S_Q$ contains elements of degree $1$ (given by
the last coordinate), we can assume that $\deg v_i=0$ for $i\in[1,m]$. In
sufficiently high degree we can find a lattice point $x$ in $\bar S_Q$ such that
$xv_i,\,xv_i^{-1} \in \bar S_P$. We have the relation $(xv_i)(xv_i^{-1})=x^2$.

It follows that $h(xv_i)=a_ix_i$, equivalently $h(x(a_i^{-1}v_i))=x_i$, for some
$x_i\in \bar S_Q$ and $a_i\in k^*$. After a toric `correction' leaving $k[\bar
S_Q]$ fixed we can assume $a_i=1$ for all $i$.

After the inversion of the elements of $S_P$, we can further extend
the homomorphism $h$ to a map defined on the Laurent polynomial ring $k[L]$.
Then we have
$$
h(v_i x x_i^{-1}) = 1.
$$
The vectors $v_i+x-x_i$ are also a basis of a complement of $U$, and thus part
of a basis of $L$. Therefore the elements
$$
v_i x x_i^{-1} - 1,\qquad  i=[1,m],
$$
generate a prime ideal of height $m$ in $k[L]$.

It is now clear that $h$ (after the toric correction) is just the retraction
$\rho_{(P,H,W)}$ where $H$ is the affine hull of $Q$ in $\RR^n$ and $W$ is the
sublattice of $\ZZ^n$ generated by the vectors $v_i+x-x_i$ upon the
identification of $\ZZ^n$ with the degree $0$ sublattice of $\ZZ^{n+1}$.
\end{proof}

Example \ref{4.2} shows that even a based retraction $h$ of $k[P]$ need not be
tame if the base does not intersect the interior of $P$ and $h$ has codimension
$\ge 2$. However, in codimension $1$ all based retraction are tame, as follows
from Theorem \ref{7.2a} and

\begin{theorem}\label{7.2b}
Suppose the codimension $1$ retraction $h:k[S_P]\to k[S_P]$ has a base $F$ not
intersecting the interior of $P$. Then $F$ is a facet of $P$ and
$h^{\epsilon}=\iota\circ\pi_F$ for some $\epsilon\in\AA(F)$ and a
$k$-algebra embedding $\iota:k[S_F]\to k[S_P]$.
\end{theorem}

In the proof we will use a general fact on pyramids.
Recall that a {\em pyramid} $\Pi\subset\RR^n$ is a polytope which is
spanned by a point $v$ and a polytope $B$ such that the affine hull of
$B$ does not contain $v$. In this situation $v$ is called an {\em apex}
and $B$ is called a {\em base} of $\Pi$.

\begin{lemma}\label{7.4}
Let $\Pi\subset\RR^n$ be a pyramid and $\Pi=\Pi_1+\Pi_2$ be a Minkowski
sum representation by polytopes $\Pi_1,\Pi_2\subset\RR^n$. Then
both $\Pi_1$ and $\Pi_2$ are homothetic images of $\Pi$ (with respect to
appropriate centers and non-negative factors).
\end{lemma}

\begin{proof}
The case $\dim(\Pi)=2$ is an easy exercise.

Now we use induction on $\dim(\Pi)$. Assume $\dim(\Pi)=n$ and assume the claim
has been shown for pyramids of dimension $\dim(\Pi)-1$. Consider any
$(n-1)$-dimensional subspace $\Lambda\subset\RR^n$ perpendicular to the
base $B\subset\Pi$. For a polytope
$R\subset\RR^n$ let $R_{\Lambda}$ denote the image of $R$ in $\Lambda$
under the orthogonal projection $\RR^n\to\Lambda$. Then $\Pi_\Lambda$ is
an $(n-1)$-dimensional pyramid and we have the Minkowski sum
representation
$$
\Pi_\Lambda=(\Pi_1)_\Lambda+(\Pi_2)_\Lambda.
$$
By induction hypothesis there are homothetic transformations of
$\Lambda$ transforming $\Pi_\Lambda$ into $(\Pi_1)_\Lambda$ and
$(\Pi_2)_\Lambda$ respectively. Considering all the possible subspaces
$\Lambda\subset\RR^n$ we conclude that
\begin{itemize}
\item[(i)]
both $\Pi_1$ and $\Pi_2$ are $n$-pyramids (provided none of them is
just a point -- in this situation the lemma is obvious) such
that the cones they span at corresponding vertices are parallel
shifts of the cone spanned by $\Pi$ at its apex $v$,
\item[(ii)]
the corresponding bases of $\Pi_1$ and $\Pi_2$ are parallel to $B$.
\end{itemize}
That is exactly what we wanted to show.
\end{proof}

\begin{proof}[Proof of Theorem \ref{7.2b}]
As in the proof of 7.2 we can assume $k[F]=\Im(h)$, and, furthermore, 
$S_P\cap\Ker(h)=\emptyset$, for otherwise $h$ itself passes through a facet 
retraction. Thus $h$ can be extended to the Laurent polynomial ring $k[L]$, 
$L=\ZZ^{n+1}$, and in particular to a retraction of $k[\bar S_P]$ with image 
$k[\bar S_Q]$. The latter restricts to retractions $k[iP]\to k[iQ]$ for all $i$. 
The kernel of the extension $h'$ is a height $1$ prime ideal and thus principal; 
$\Ker(h')=\phi K[L]$ and $\Ker(h)=(\phi k[L])\cap k[P]$ for some element 
$\phi\in k[L]$.

Since $F$ is a base of $h$, $\Ker(h)$ contains the elements $x-\ell$,
$x\in\L_P\setminus F$, $\ell=h(x)$, and $\ell$ is a linear form on the points of $\L_F$. Then
$\N(\phi)$ is a Minkowski summand of the pyramid $\N(x-\ell)$ with vertex at
$x$. One can shift $\N(\phi)$ by an integer vector
into $\N(x-h(x))\subset P$ such that the image $R$ satisfies
\begin{equation}
\tag{$**$}R\subset P\ \text{and}\ R\cap F\neq\emptyset.
\end{equation}
Evidently $R$ is the Newton polytope of $y\phi$ for some $y\in\ZZ^{n+1}$.
Replacing $\phi$ by $y\phi$, we can assume that $\N(\phi)$ satisfies $(**)$.

By Lemma \ref{7.4} $\N(\phi)$ is homothetic to $\N(x-\ell)$. Clearly,
$F\cap\N(\phi)$ is a base of $\N(\phi)$. The corresponding apex of
$\N(\phi)$ is some $z\in\L_P\setminus F$.

Consider the valuation
$$
v_F:\ZZ^{n+1}\to\ZZ
$$
determined by the conditions:
$$
\Im(v_F)=\ZZ,\ v_F(\L_F)=0,\ v_F(\L_P)\geq0
$$
We claim: $v_F(z)=1$ and $y+b-z\in P$ for any $b\in\L_{F\cap\N(\phi)}$
and $y\in\L_P\setminus F$.

In fact, for $i\in\NN$ big enough there is an element $z'\in\L_{iP}$ such
that $v_F(z')=1$. Since $iF$ is a base of the induced retraction $\bar
h_i:k[{iP}]\to k[{iP}]$ there exists a linear form $\ell'$ on $\L_{iF}$ such that $z'-\ell'\in\Ker(h)$.
Thus $\N(\phi)$ is a Minkowski summand of $\N(z'-\ell')$. Because
of the condition $(**)$ we conclude $v_F(z)\leq v_F(z')$. Hence
$v_F(z)=1$.

Now choose $y\in\L_P\setminus F$. Since $F$ is a base of $h$, we can write
$y-\ell''\in\Ker(h)$ for some linear form $\ell''$ on the points of $\L_F$.
Therefore, the pyramid $\N(\phi)$ is a Minkowski summand of the pyramid
$\N(y-\ell'')$ which has its apex at $y$. By Lemma \ref{7.4} the cones
spanned by these pyramids at their vertices are the same modulo a parallel
shift. This observation in conjunction with the already established equality
$v_F(z)=1$ makes the claim clear.

We have shown that the vectors $b-z\in\ZZ^n$, $b\in\L_{F\cap\N(\phi)}$,
are column vectors for $P$. Now, by Lemma \ref{6.1}(b) there exists
$\epsilon\in\AA(F))$ such that $\epsilon(\phi)=cz$ for some $c\in
k^*$. Therefore, $\Ker(h^{\epsilon})$ is the monomial prime ideal
$(\L_P\setminus F)k[P]\subset k[P]$, and this finishes the proof of
Theorem \ref{7.2b}.
\end{proof}

Using Theorem \ref{7.2a} and Theorem \ref{7.2b} and Borel's theorem on maximal
tori in linear groups \cite{Bor} we now derive another sufficient condition for a
codimension 1 retraction to be tame.

\begin{theorem}\label{7.5}
A codimension $1$ retraction $h:k[P]\to k[P]$ is tame
if the graded automorphism group of the $k$-algebra $\Im(h)$ contains a
closed $\dim(P)$-torus whose action on $\Im(h)$ extends to an algebraic
action on $k[P]$.
\end{theorem}

\begin{proof}
Let $\TT\subset\ggr.aut(\Im(h))$ be such a closed torus.
Then there is an algebraic embedding $\TT\subset\Gamma_k(P)$ and a maximal
intermediate torus $\TT\subset\TT'\subset\Gamma_k(P)$. By Borel's theorem
$$
\TT'=\TT_k(P)^\alpha=\{\alpha\circ\tau\circ\alpha^{-1}|\ \tau\in\TT_k(P)\}
$$
for some $\alpha\in\Gamma_k(P)$. Let us
show that $h^\alpha$ is a based retraction. By the previous theorems this
completes the proof. We can assume that
$\Ker(h^\alpha)\subset k[P]$ is not a monomial ideal, or, equivalently,
$\Ker(h^\alpha)$ contains no monomials. Otherwise $h$ factors through a facet projection, and we are done.
\medskip

\noindent{\em Claim.} $\Ker(h^\alpha)\subset k[P]$ is a binomial ideal.
\medskip

We know that $\Ker(h^\alpha)$ is stabilized by the codimension 1 subtorus
$\TT^\alpha\subset\TT_k(P)$. Since $\Ker(h^\alpha)$ a is non-monomial
height 1 prime ideal of $k[P]$ there is an element $\phi\in k[\gp(S_P)]$
such that
$$
\Ker(h^\alpha)=k[P]\cap(\phi k[\gp(S_P)]).
$$
(Compare with the proof of Theorem \ref{5.1}.) In this situation for any
$\tau\in\TT^\alpha$ we have $\tau(\phi) k[\gp(S_P)]=\phi k[\gp(S_P)]$.
Since $\TT^\alpha$ is a subtorus of the embedded torus, the latter
equation is equivalent to the existence of an element $a_\tau\in k^*$
such that $\tau(\phi)=a_\tau\phi$. ($\TT^\alpha$ acts naturally on
coordinate ring $k[\TT_k(P)]=k[\gp(S_P)]$.) This
means that for any $\tau\in\TT^\alpha$ the terms in the canonical
$k$-linear expansion of $\phi$ are multiplied by the same scalar from
$k^*$ when $\tau$ is applied, a condition equivalent to
to the segmentonomiality of $\phi$. In fact, if $\dim \N(\phi)\geq2$
we would have three non-collinear terms $T_1$, $T_2$ and $T_3$ in the
$k$-linear expansion of $\phi$, and the condition
$$
\frac{\tau(T_1)}{T_1}=\frac{\tau(T_2)}{T_2}=\frac{\tau(T_3)}{T_3}
$$
on $\tau\in\TT_k(P)$ would be equivalent to the condition that each $\tau$ is
a solution to {\em two independent} binomial equations on
$\TT_k(P)=(k^*)^{\dim(P)+1}$. Hence, $\TT^\alpha\subset\TT_k(P)$ would
be a sub-torus of codimension at least $2$ -- a
contradiction. By this argument we have shown that $\phi$ is a
segmentonomial and, hence, a binomial by Theorem \ref{5.1}(a).
Therefore, by Theorem 3.4 we can assume $\Im(h)=k[S]$, $S=h(S_P)$. In
particular,
$$
\sigma_*(\TT^\alpha)=\TT_k(S),
$$
where $\sigma_*:\ggr.aut(\Im(h))\to\ggr.aut(k[S])$ is the natural map
induced by $\sigma$ and $\TT_k(S)$ is the embedded torus of
$\Spec(k[S])$ corresponding to $S$. The equality holds true because
monomial preserving automorphisms are mapped to monomial preserving automorphisms.

Identifying $\Im(h)$ and $k[S]$ via $\sigma$ and $\TT^\alpha$ and
$\TT_k(S)$ via $\sigma_*$ we get a $\TT^\alpha$-equivariant embedding
$$
\iota:k[S]\to k[P].
$$
Assume $\{s_1,\ldots,s_m\}$ be the set of degree 1 terms in
$k[S]$. Then for any $\tau\in\TT^\alpha$ there exist $a_\tau^{(i)}\in
k^*$ such that
$$
\tau(\iota(s_i))=a_\tau^{(i)}\iota(s_i),\qquad i\in[1,m].
$$
Arguments similar to those in the proof of the claim above show that the
polytopes $\N(\iota(s_i))$, $i\in[1,m]$, are segments parallel to $\N(\phi)$,
maybe some of them degenerated into points.

Let $z\in\RR\otimes\gp(S_P)\setminus\{0\}$ be such that the line $\RR
z\subset\RR\otimes\gp(S_P)$ is parallel to all the $\N(\iota(s_i))$.
Denote by $t_i$, $i\in[1,m]$ the upper end-points of $\N(\iota(s_i))$ in
the direction $\RR_+z$. Then standard arguments with Newton polytopes
ensure that the $t_i$ are subject to the same affine {\em binomial}
dependencies as the $s_i$. It is now clear that the
subpolytope
$$
\conv(\{t_1,\ldots,t_m\})\subset P
$$
is a base for the retraction $h:k[P]\to k[S]$.
\end{proof}

\section{Retractions of polygonal algebras}\label{gon}

Throughout this section $P$ denotes a lattice polygon (i.~e. $\dim(P)=2$).

\begin{theorem}\label{8.1}
Any codimension $1$ retraction $h:k[P]\to k[P]$ is based and,
therefore, either $h^{\tau}=\iota\circ\rho_{(P,H,w)}$ for some lattice
segmental fibrations $(P,H,w)$, $\tau\in\TT_k(P)$ and a $k$-embedding
$\iota:k[{H\cap P}]\to k[P]$, or $h^{\epsilon}=\iota\circ\pi_F$
for a facet $F\subset P$, $\epsilon\in\EE_k(P)$ and a $k$-embedding
$\iota:k[F]\to k[P]$. In particular, $h$ is tame.
\end{theorem}

By Theorems \ref{7.2a} and \ref{7.2b} it is enough to find a base for $h$. The
first step in its construction is given by

\begin{proposition}\label{8.2}
A multiple $c\Delta_1$, $c\in\NN$, of the unit segment $\Delta_1$ can be
embedded as a lattice polytope into $P$ if and only if there is a
$k$-algebra embedding of $k[{c\Delta_1}]$ into $k[P]$.
\end{proposition}

\begin{proof}
Clearly, without loss of generality we can assume $c\geq2$.

Let $\epsilon:k[{c\Delta_1}]\to k[P]$ be an embedding. We write
$$
\L_{c\Delta_1}=\{x_0,x_1,\ldots,x_c\}.
$$
Thus we have the equations
$
\epsilon(x_{i-1})\epsilon(x_{i+1})=\epsilon(x_i)^2
$
for $i\in[1,c-1]$. Put
$$
\phi=\frac{\epsilon(x_1)}{\epsilon(x_0)}=\frac{\epsilon(x_2)}{\epsilon(x_1)}=\cdots.
$$
In the Laurent polynomial ring $k[\gp(S_P)]=k[\ZZ^3]$ we can write
$
\phi=\phi_1/\phi_2
$
with coprime $\phi_1,\phi_2$. The
equality
$
\phi_2^c\epsilon(x_c)=\phi_1^c\epsilon(x_0)
$
(and the factoriality of $k[\ZZ^3]$) imply that $\phi_1^c$ divides $\epsilon(x_c)$
and $\phi_2^c$ divides $\epsilon(x_0)$.
\medskip

\noindent{\em Case} (a). Both $\phi_1$ and $\phi_2$ are monomials in $k[\ZZ^3]$.
In this situation the Newton polygon $\N(\epsilon(x_c))$ is the parallel shift
of $\N(\epsilon(x_0))$ by the $c$-th multiple of the vector representing the
support term of the monomial $\phi$. But then the existence of the desired
embedding $c\Delta_1\to\conv(\N(\epsilon(x_0)),\N(\epsilon(x_c))\subset P$ is
obvious.
\medskip

\noindent{\em Case} (b). At least one of $\phi_1$ and $\phi_2$, say $\phi_1$, is
not a monomial. Then $c\Delta_1$ can be embedded in any of the edges of the
polygon $\N(\phi_1^c)=c\N(\phi_1)$.
Since $\epsilon(x_c)=\psi\phi_1^c$ for some $\psi\in k[\ZZ^3]$, we get
$$
\N(\epsilon(x_c))=\N(\psi)+\N(\phi_1^c)
$$
and the existence of an embedding
$c\Delta_1\to\N(\epsilon(x_c))\subset P$ is evident.
\end{proof}

\begin{remark}\label{8.3}
We expect that Proposition \ref{8.2} holds without any restrictions: for lattice
polytopes $P$ and $Q$ a $k$-algebra embedding $k[Q]\to k[P]$ should only exist
if $Q$ can be embedded into $P$ (as a lattice subpolytope).
\end{remark}

Below we will need the following notion. For a lattice polygon
$P\subset\RR^2$ its {\em lattice width} $\width_\ell(P)$ with respect to the
line $\ell=z'+\RR z$, $z,z'\in\ZZ^2$, $z\neq0$, is defined as
$$
\width_\ell(P)=\max_{x\in P\cap\ZZ^2}\phi(x)-\min_{y\in P\cap\ZZ^2}\phi(y)
$$
where $\phi:\ZZ^2\to\ZZ$ is a $\ZZ$-linear form with $\phi(z)=0$ and $\phi(\ZZ)=\ZZ$. In other words, $\width_\ell(P)$ is the number of integral translates of $\ell$
intersecting $P$ minus $1$.

\begin{proof}[Proof of Theorem \ref{8.1}.] Let $h:k[P]\to k[P]$ be a codimension 1
retraction. By Theorem \ref{2.2} there is an isomorphism $\Im(h)\approx k[{c\Delta_1}]$ for some
$c\in\NN$, and Proposition \ref{8.2} then yields an embedding $c\Delta_1\to
P$. Let $Q$ denote its image.

If the restriction of $h$ to the subalgebra $k[Q]\subset k[P]$ is
injective then it is evidently bijective, and we are done by Lemma \ref{7.1}.
Therefore we can assume
$$
\Ker(h)\cap k[Q]\neq\{0\}.
$$
Clearly, $I=\Ker(h)\cap k[Q]\subset k[Q]$ is segmentonomial. Since
$\het(\Ker(h))=1$ the ideal $\Ker(h)$ is a minimal prime ideal of $I$. If
$\Ker(h)$ contains a monomial, then $h$ factors through a facet projection
$\pi_F$, and then $F$ is a base for $h$. In the other case it is enough to show
there is yet {\em another} embedding $c\Delta_1\to P$, say with image $R$, such
that $Q$ and $R$ are not parallel. Theorem \ref{5.1}(b) shows that we then have
found a base.

Let $\ell$ denote the line spanned by $Q$. By Theorem \ref{5.1}(b) we have
\begin{equation}
\tag{$*$}\width_\ell(P)=c.
\end{equation}
We use the notation
$$
P_i=\N(\epsilon(x_i)),\qquad i\in[0,c].
$$
(As in the proof of Proposition \ref{8.2} we write
$\L_{c\Delta_1}=\{x_0,x_1,\ldots,x_c\}$.) Evidently there is no loss in
generality in assuming that
$$
P=\conv(P_0\cup\dots\cup P_c)
$$
since, roughly speaking, the retraction can be restricted to the lattice
polytope on the right hand side.\smallskip

\noindent{\em Case} (a): $P_0 \neq P_1$.
For every linear form $\lambda:\RR^2\to\RR$ we consider all triples
$(\lambda,v_0,v_1)$ such that $v_0\in P_0$, $v_1\in P_1$ and
$$
\lambda(v_i)>\lambda(x)\qquad\text{for all }x\in P_i,\quad x\neq v_i,\quad i=0,1.
$$
In particular $v_0$ and $v_1$ are vertices of $P_0$ and $P_1$ respectively. We have
$$
iP_1=(i-1)P_0+P_i,\qquad i\in[2,c].
$$
It follows that $\lambda$ has a unique maximum on $P_i$, necessarily taken at
$v_i=iv_1-(i-1)v_0$, $i\in[0,c]$. Thus we get the system
$$
\{v_0,v_1,\ldots,v_c\}\subset P
$$
of collinear lattice points.

\noindent{\em Subcase} $(\mathrm a_1)$. There is a triple $(\lambda,v_0,v_1)$
such that $[v_0,v_1]$ is {\em not parallel to $\ell$}, and in particular
$v_0\neq v_1$. Then the points $v_0,v_1,\ldots,v_c$ are pairwise different, and
so we find $c+1$ successive lattice points in $P$ on a line not parallel to
$\ell$.

\noindent{\em Subcase} $(\mathrm a_2)$. In the other case one notes first that
the edges of $P_0$ and $P_1$ are pairwise parallel under an orientation
preserving correspondence (in other words, $P_0$ and $P_1$ have the same normal
fans; see \cite{GKZ} for this notion). Moreover, if $v_0$ is a vertex of $P_0$
and $v_1$ is the corresponding vertex of $P_1$, then $v_0=v_1$ or the segment
$[v_0,v_1]$ is parallel to $\ell$.

Let us choose $z\in\RR^2$, $z\neq 0$ such that $\ell$ is parallel to $\RR z$,
and let $\xi$ be a non-zero linear form with $\xi(z)=0$.

Furthermore we consider the {\em upper} and {\em lower boundary} of $P_0$ in
direction $z$. (A point $y$ on $\partial P_0$ belongs to the lower boundary if
the ray in direction $z$ through $y$ enters $P$ at $y$, and the upper boundary
is defined analogously; see Figure \ref{boundaries}.) Each of these subsets has
two endpoints $u_0^-, u_0^+$ and $l_0^-,l_0^+$ respectively (which may
coincide). We choose the indices such that $\xi(u_0^-)\le \xi(u_0^+)$ and
$\xi(l_0^-)\le \xi(l_0^+)$. It is clear that $\xi$ takes its minimum on $P_0$ in
$u_0^-$ and $l_0^-$ and its maximum in the other two points.
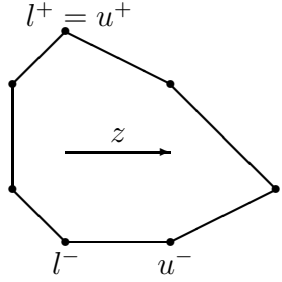
\begin{figure}[htb]
\begin{center}
\begin{picture}(6.5,4)(0,1)
\put (0,2){\vertex}
\put (0,4){\vertex}
\put (1,5){\vertex}
\put (3,4){\vertex}
\put (5,2){\vertex}
\put (3,1){\vertex}
\put (1,1){\vertex}
\put (3,1){\line(-1,0){2}}
\put (1.0,2.7){\vector(1,0){2}}
\put (1.85,2.9){$z$}
\put (0.75,0.4){$l^-$}
\put (2.75,0.4){$u^-$}
\put (0.25,5.1){$l^+=u^+$}
\thicklines
\put (0,4){\line(1,1){1}}
\put (1,5){\line(2,-1){2}}
\put (3,4){\line(1,-1){2}}
\put (5,2){\line(-2,-1){2}}
\put (1,1){\line(-1,1){1}}
\put (0,2){\line(0,1){2}}
\end{picture}
\end{center}
\caption{Lower and upper boundary}\label{boundaries}
\end{figure}
Applying the same construction to $P_1$ we obtain the points $u_1^-,u_1^+$ and
$l_1^-,l_1^+$. The correspondence between the vertices of $P_0$ and those of
$P_1$ pairs $u_0^-$ with $u_1^-$ etc.

If $u_0^-=u_1^-$ and $l_0^-=l_1^-$, then $P_0=P_1$ by our initial
observation in this subcase. But we are assuming that $P_0\neq P_1$ in Case (a).

For example, suppose that $u_0^-\neq u_1^-$. Then we can find a linear form
$\lambda$ such that $u_0^-$ and $u_1^-$ are the unique maxima of $\lambda$ on
$P_0$ and $P_1$ respectively. So the triple $(\lambda,u_0^-,u_1^-)$ defines
$c+1$ lattice points $v_0,\dots,v_c$ in $P$ on a line $\ell'$ parallel to
$\ell$, and one sees easily that $\ell'$ intersects $P$ in an edge: it is
impossible that $\xi(t)< \xi(u_0^-)=\xi(u_1^-)$ for some $t\in P_i$, $i\in
[0,c]$.

Since
$\width_\ell(P)=\width_{\ell'}(P)=c$ there is a lattice point $t\in P$ such that
the triangle $T\subset P$ spanned by $t$ and $c+1$ successive lattice points
$$
w_0,w_1,\ldots,w_c\in\ell'\cap \conv(\{v_0,v_1,\ldots,v_c\})\subset P
$$
has $\width_\ell(T)=c$. We apply an
integral affine transformation $\psi:\RR^2\to\RR^2$ such that
$\psi(t)=(0,0)$, and
$$
\psi(\conv(\{w_0,v_1,\ldots,w_c\}))=\conv((c,r),(c,r+c))
$$
for some $r\in\ZZ$. There exists $d\in\ZZ$ such that $dc\in [r,r+c]$. The
segment $[(0,0),(c,dc)]$ clearly contains exactly $c+1$ equidistant lattice
points and is not parallel to the image of $\ell$ under the transformation
$\psi$. Again we are done.
\medskip

\noindent{\em Case} (b): $P_0=P_1$. Then it follows from the equation
$iP_1=(i-1)P_0+P_i$, $ i\in[2,c]$ that all the $P_i$ coincide and are therefore
identical with $P$. If $P$ has an edge containing at least $c+1$ lattice points,
then either this edge is not parallel to $Q$ or there is a triangle $T$ as in
subcase $(\mathrm a_2)$, and in both cases we are done.

Only the case in which every edge of $P$ has lattice length strictly less than
$c$ is left. In this situation we resort to the general result below
(Proposition \ref{8.4}) which we have singled out because of its independent
interest.
\end{proof}

\begin{proposition}\label{8.4}
Assume the edges of the lattice polygon $P\subset\RR^2$ have lattice
lengths strictly less than $c$ and $\width_\ell(P)\leq c$ for some line
$\ell=\RR x\subset\RR^2$, $x\in\ZZ^2$, $x\neq0$. Then there is
no $k$-algebra embedding $\epsilon:k[{c\Delta_1}]\to k[P]$ such that
$$
P=\N(\epsilon(x_0))=\N(\epsilon(x_1))=\cdots=\N(\epsilon(x_c))
$$
where $x_0,\dots,x_c$ denote the successive lattice points of $c\Delta_1$.
\end{proposition}

\begin{proof}
Assume to the contrary that such $\epsilon$ exists. For a subset $W\subset P$
and an element
$$
f=a_1y_1+\cdots+a_ry_r\in k[P],\ \ a_j\in k,\ \ y_j\in\L_P,
$$
we put
$$
f_W=a_{j_1}y_{j_1}+\cdots+a_{j_t}y_{j_t}, \qquad
\{y_{j_1},\ldots,y_{j_t}\}=W\cap\{y_1,\ldots,y_r\}.
$$
Let $E\subset P$ be an edge. It is impossible to find a graded
$k$-algebra embedding of $k[{c\Delta_1}]$ into
$k[E]$ (by Hilbert function reasons or Proposition \ref{8.2}). In particular,
the homomorphism
$$
\pi_E\circ\epsilon:k[{c\Delta_1}]\to k[E],\qquad x_i\mapsto\epsilon(x_i)_E,
$$
is not injective. (As usual, $\pi_E$ is the facet projection $k[P]\to
k[E]$.)
It follows that
$$
(\pi_E\circ\epsilon)(k[{c\Delta_1}])\approx k[\ZZ_+].
$$
as graded $k$-algebras. Thus there is an element $\xi\in k^*$ such that
\begin{equation}
\tag{1}\epsilon(x_0)_E=\xi\epsilon(x_1)_E=\xi^2\epsilon(x_2)_E=\cdots=\xi^c\epsilon(x_c)_E.
\end{equation}
Consider the toric automorphism $\tau:k[{c\Delta_1}]\to
k[{c\Delta_1}]$ determined by
$$
x_i\mapsto\xi^i x_i,\qquad i\in[0,c].
$$
Then
$$
\epsilon'(x_0)_E=\epsilon'(x_1)_E=\cdots=\epsilon'(x_c)_E,
$$
where $\epsilon'=\epsilon\circ\tau$.

Applying the same arguments to the embedding $\epsilon':k[{c\Delta_1}]\to k[P]$
and an edge $E'\subset P$ sharing an end-point with $E$ we conclude
$$
\epsilon'(x_0)_{E'}=\epsilon'(x_1)_{E'}=\cdots=\epsilon'(x_c)_{E'}.
$$
(Because of the common end-point there is no need to further `correct'
$\epsilon'$ by a toric automorphism.) We can assume $\epsilon=\epsilon'$.

Going once around $P$ we finally arrive at the equations
\begin{equation}
\tag{2}\epsilon(x_0)_{\partial P}=\epsilon(x_1)_{\partial P}=\cdots=
\epsilon(x_c)_{\partial P},
\end{equation}
where $\partial P$ denotes the boundary of $P$.

After an integral affine change of coordinates we may assume $\ell=\RR(0,1)$ and
that $P\subset\RR^2_+$. Then $k[P]$ can be considered as a subalgebra of
$k[X,Y,Z]$ generated by monomials $X^iY^jZ$. Moreover, after a parallel
translation of $P$ we can assume that $P$ contains a lattice point $(0,j)$, or
equivalently, $k[P]$ contains a monomial $Y^jZ$.

For any element $a\in k^*$ we have the retraction
$$
\pi^a:k[X,Y,Z]\to k[X,Z],
\qquad
\pi_a(X^iY^jZ^m)\mapsto a^jX^iZ^m
$$
for all $i,j,m$. The image of $k[P]$ under $\pi^a$ is a polytopal algebra
$k[c'\Delta_1]$ generated by the elements $Z,XZ,\dots,X^{c'}Z$ with
$c'=\width_\ell(P)$. 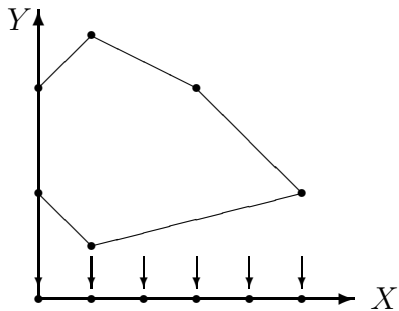
\begin{figure}[htb]
\begin{center}
\begin{picture}(6.5,5)(0,0)
\put (0,2){\vertex}
\put (0,4){\vertex}
\put (1,5){\vertex}
\put (3,4){\vertex}
\put (5,2){\vertex}
\put (1,1){\vertex}
\put (0,4){\line(1,1){1}}
\put (1,5){\line(2,-1){2}}
\put (3,4){\line(1,-1){2}}
\put (5,2){\line(-4,-1){4}}
\put (1,1){\line(-1,1){1}}
\multiput (0,0.8)(1,0){6}{\vector(0,-1){0.6}}
\multiput (0,0)(1,0){6}{\vertex}
\put (-0.6,5.1){$Y$}
\put (6.3,-0.2){$X$}
\thicklines
\put (0,0){\vector(0,1){5.5}}
\put (0,0){\vector(1,0){6}}
\end{picture}
\end{center}
\caption{The map $\pi^a$}\label{Last}
\end{figure}
Due to the infinity of the field $k$ there is $a\in k^*$ such that the
elements
$$
\pi^a(\epsilon(x_i))\qquad i\in[0,c],
$$
are pairwise different. In fact, it is clear that for any finite set of
polynomials ${\mathcal P}\subset k[X,Y,Z]$ there is an element $a\in k^*$ such
that $\pi^a$ maps $\mathcal P$ to a set of $\#\mathcal P$ pairwise different
polynomials from $k[X,Y]$.

Set
$$
\pi^a(\epsilon(x_j))=\sum_{i=0}^{c'} \alpha_{ij} X^iZ.
$$
In view of the fact that the elements $\epsilon(x_j)$ coincide `along $\partial
P$', we see that the coefficients $\alpha_{0j}$ coincide for all $j$, and the
same holds for the coefficients $\alpha_{c'j}$. The $(c'+1)\times (c+1)$ matrix
$(\alpha_{ij})$ has two linear dependent rows, and since $c'\le c$, its rank is
$<c+1$.

It follows that $\pi^a\circ\epsilon$ is not injective. Thus we see as above that
there exists $\eta$ in $k^*$ with
$$
\pi^a(\epsilon(x_0))=\eta\pi^a(\epsilon(x_1))=\eta^2\pi^a(\epsilon(x_2))=\cdots=
\eta^c\pi^a(\epsilon(x_c)).
$$
Then however $\alpha_{00}=\eta\alpha_{01}=\dots=\eta^c\alpha_{0c}$. By a
suitable choice of $a$ we can furthermore achieve that $\alpha_{00}\neq0$. Then
it follows immediately that $\eta=1$. This is a contradiction since we have
chosen $a$ such that the $\pi^a(\epsilon(x_i))$ are pairwise different.
\end{proof}

\begin{remark}\label{8.6}
Essentially the same arguments as we have used in \S7 and \S8 imply the
analogous results for graded retractions of the semigroup algebras corresponding
to affine positive {\em normal} semigroups (we only assume all monomials are
homogeneous of positive degree) -- first one observes that there is the
appropriate notion of a lattice fibration for rational polyhedral pointed cones,
and then the arguments are applied to all Veronese subalgebras. (They are
essentially the same as polytopal semigroup rings.)

This remark is analogous to Remark 3.3(c) in \cite{BG1}, saying that the
description of graded automorphism groups for graded affine normal
semigroup rings is essentially the same as that for polytopal semigroup
rings.
\end{remark}

\end{document}